\documentclass[12pt]{article}
\usepackage{setspace, graphicx, color,epsfig,graphics,amssymb,amsmath,subfigure,verbatim,lscape,fullpage}
\usepackage[all]{xy}
\pdfoutput=1

\ifx\pdfoutput\undefined
	\usepackage[hypertex]{hyperref}
\else
	\usepackage[pdftex,hypertexnames=false]{hyperref}
\fi

    \def \R {{\mathbb{R}}}
   \def \C {{\mathbb{C}}}
   \def \N {{\mathbb{N}}}
   \def \Z {{\mathbb{Z}}}
    \def \e {{\varepsilon}}

    \def \o {{\theta}}

    \def \y {{y_1}}

    \def \D {{\Delta}}
    \def \d {{\delta}}
    
    \def \w {{w_{\tilde{\y}}}}

    \def \l {{\lambda}}

    \def \n {{\eta}}
    \def \la {{\langle}}
    \def \ra {{\rangle}}

\newcommand{\beq}{\begin{equation}}
\newcommand{\eeq}{\end{equation}}
\newcommand{\beqq}{\begin{equation*}}
\newcommand{\eeqq}{\end{equation*}}
\newcommand{\beqr}{\begin{eqnarray}}
\newcommand{\eeqr}{\end{eqnarray}}
\newcommand{\beqrn}{\begin{eqnarray*}}
\newcommand{\eeqrn}{\end{eqnarray*}}
\newcommand{\beqn}{\begin{equation*}}
\newcommand{\eeqn}{\end{equation*}}
\newcommand{\bei}{\begin{itemize}}
\newcommand{\beii}{\begin{itemize} \item}
\newcommand{\eei}{\end{itemize}}
\newcommand{\bmei}{\begin{itemize} \compactlist}
\newcommand{\emei}{\end{itemize}}
\newcommand{\ben}{\begin{enumerate}}
\newcommand{\een}{\end{enumerate}}
\newcommand{\bes}{\begin{small}}
\newcommand{\ees}{\end{small}}
\newcommand{\bec}{\begin{center}}
\newcommand{\eec}{\end{center}}

\newcommand{\eps}{\epsilon}

\newtheorem{lemma}{Lemma}

\begin{document}
\title{Shared inputs, entrainment, and desynchrony in elliptic bursters:  from slow passage to discontinuous circle maps}

\author{Guillaume Lajoie and Eric Shea-Brown}
\date{\today}
\maketitle

\begin{abstract}  What input signals will lead to synchrony vs. desynchrony in a group of biological oscillators?  This question connects with both classical dynamical systems analyses of entrainment and phase locking and with emerging studies of stimulation patterns for controlling neural network activity.  Here, we focus on the response of a population of uncoupled, elliptically bursting neurons to a common pulsatile input.  We extend a phase reduction from the literature to capture inputs of varied strength, leading to a circle map with discontinuities of various orders.  In a combined analytical and numerical approach, we apply our results to both a  normal form model for elliptic bursting and to a biophysically-based neuron model from the basal ganglia.  We find that, depending on the period and amplitude of inputs, the response can either appear chaotic (with provably positive Lyaponov exponent for the associated circle maps), or periodic with a broad range of phase-locked periods.  Throughout, we discuss the critical underlying mechanisms, including slow-passage effects through Hopf bifurcation, the role and origin of discontinuities, and the impact of noise.  
\end{abstract}


\section{Introduction}

Many types of physical and biological systems exhibit intrinsic bursting -- rapid discharges of consecutive, fast dynamical events separated by periods of quiescence.  In particular, bursting neurons serve myriad functions in the nervous system; prominent among these is their role in central pattern generators that create rhythmic neural activity \cite{CalinJageman:2007p7271,Cymbalyuk:2002p7245,shilnikov:037120,Coombes:2005p3301}.  Bursting dynamics also feature in pathological oscillations associated with disease conditions, as for basal ganglia networks in Parkinsons disease~\cite{Pirini:2009p521,Terman:2002p125,Best:2007p128,Azad:2010p4782,Izhikevich:2001p2437}, where elevated synchrony and rhythmicity among neurons is linked to motor symptoms.  

Here, we focus on synchrony and desynchrony among bursting neurons in the simplest possible setting:  a population of uncoupled bursting neurons receiving a common input signal.  We study elliptic bursters \cite{Izhikevich:2000p2513,Rinzel:1989p2766} -- that is, non-linear oscillators with fast and slow variables, and for which burst onset is caused by passage through a (subcritical) Hopf bifurcation in the fast subsystem and burst offset follows from a saddle-node bifurcation of limit cycles (see Sect. \ref{geometry} below).  The driving signals are periodic pulsatile inputs.  

\begin{figure}[h]
\begin{center}
\includegraphics{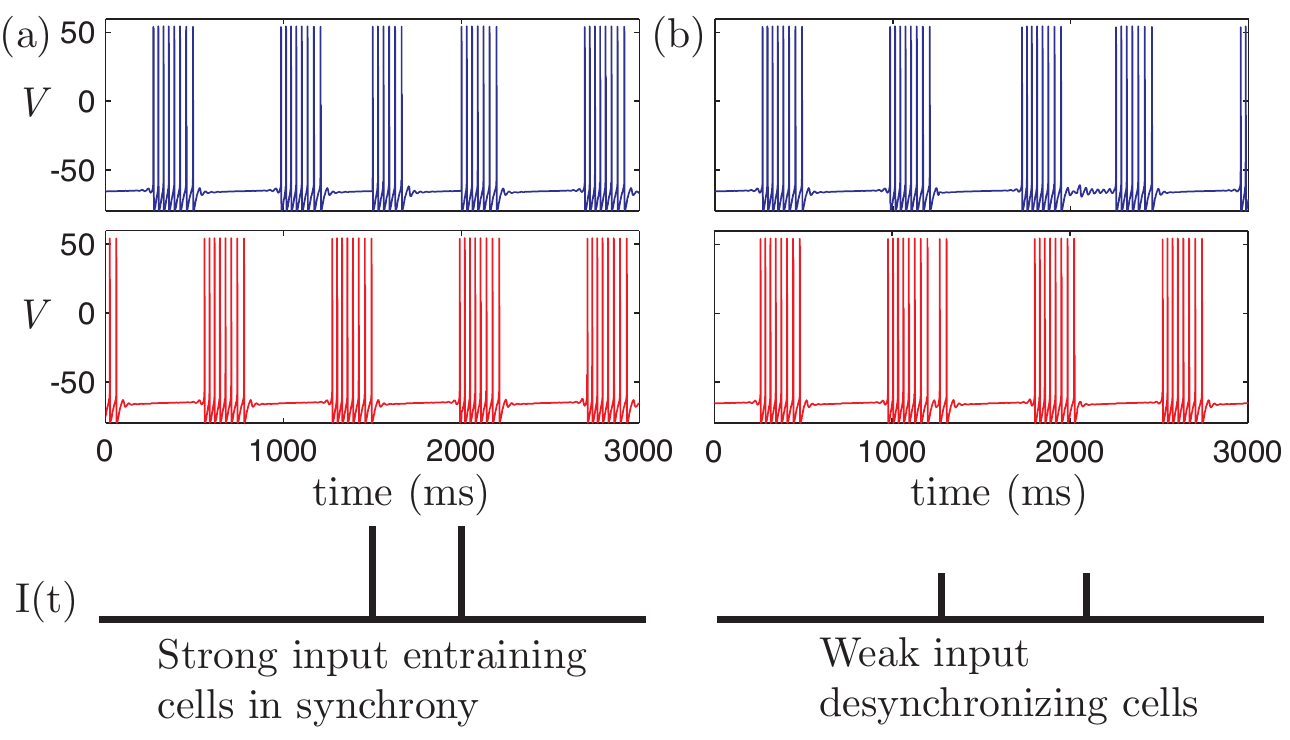}
\caption{Simulated voltage traces of two uncoupled cells receiving common inputs (from the conductance-based neuron model of Eqn.~(\ref{eqn:Terman}); see Sect.~\ref{numerics}).  (a) A strong input synchronizes cells that are initially out of phase. (b) A weak input desynchronizes cells with close by initial phases.}
\label{fig:synch_desynch}
\end{center}
\end{figure}

We find that there is rich variety in the response to these inputs, depending on their strength and frequency. An illustrative example is presented in Fig.~\ref{fig:synch_desynch}, where we plot simulated voltage traces of two bursting cells, both receiving a common pulsatile input $I(t)$. In the left panel, the cells' bursting phases are initially well separated, but a pair of ``strong" input pulses synchronizes them. In the right panel, the cells' phases are initially nearly synchronized, but a relatively ``weak" input drives them apart.  As we will show, the outcome depends on pulse strength and on inter-pulse timing in interesting ways that arise directly from the dynamics of elliptic bursting.

Our goal is to understand the mechanisms responsible for this and related phenomena. 
Specifically, we ask whether or not cells will entrain or become separated under the driving effect of a common periodic input signal.  Our main goal is to develop and explain a general answer to this question. We note two broad scientific motivations for doing so, although exploring the implications of our results for each largely remains for future work.  First, entrainment of a population of uncoupled cells to an input signal determines the ``reliability" of a neuron's response -- that is, the repeatability of a response  to a fixed input signal on multiple trials in which the neuron is in a different initial condition.  Reliability is fundamental in understanding how neural dynamics encode, e.g., sensory signals~\cite{Bry+76,MS95,HunterMTC98,01,LinSY07a,Banerjee:2008p6505}.  Second, common input signals to populations of bursting cells (of elliptic and other types) occur naturally in layered neural networks, as in the basal ganglia; in this brain area, common, pulsatile electrical stimuli are also artificially applied as a therapy for Parkinson's disease~\cite{Pirini:2009p521,Terman:2002p125,Best:2007p128}.  We give a very brief application to this general setting in Sect.~\ref{composition}.   

Throughout this paper, we use a normal form model for elliptic bursting developed by Izhikevich \cite{Izhikevich:2000p2513}. It is the simplest system that captures the dynamical features of elliptic bursting and we show that it accurately describes a more complex bursting neuron model derived from basal ganglia physiology.  Thus, we will often refer to the \textit{normal form model} as describing a ``cell."  

Building from, e.g.~\cite{Best:2007p128}, we show that the dynamics of the normal form model under periodic, pulsatile inputs admits an accurate reduction to a discontinuous circle map that can be analytically defined.  This map forms the basis for our theoretical results, and correctly predicts synchrony and desynchrony for the physiological model.  Nevertheless, because our results are linked to the normal form model, we note that our approach and results have potential applications in the study of general slow/fast oscillators undergoing a delayed bifurcation, well beyond neuroscience.
  
The paper proceeds as follows. Section \ref{geometry} deals with the analysis of elliptic bursting dynamics as well as phase reduction to a discrete dynamical system on the circle. Section \ref{dynamics} presents an analysis of the reduced dynamics and explores synchronizing and desynchronizing effects of common pulsatile inputs.  Next, in Section \ref{noise}, we study the effect of noise, which has a non-trivial and interesting impact on the circle map and resulting patterns of synchrony.  In Section \ref{numerics}, we carry out a series of numerical experiments that verify our reduced models. Finally, in Section \ref{composition}, we show how the circle map framework can be used to analyze the effect of multiple sequences of pulsatile inputs. As a proof of concept, we present a brief example in which an input signal is designed to compete with the effect of a global entraining drive synchronizing a population, as in the basal ganglia setting described above. 
{\it Our principal finding -- unifying analyses of several models and settings -- is that a population of (elliptic) bursting cells will desynchronize in the presence of weak to moderately strong common inputs, if these inputs have a frequency slightly slower than the natural burst frequency.}  


\section{Geometry of forced elliptic bursters}
\label{geometry}

In order to better understand the effect of a common stimulus on a population of bursters, we first describe the dynamics of single cells. Dynamical models that capture bursting usually include multiple timescales and are often called \textit{slow/fast} systems. Indeed, most intrinsically bursting solutions arise from the evolution of one or more \textit{slow} variables that periodically steer \textit{fast} variables into distinct dynamical regimes -- here, spiking and resting. 

\subsection{Timescale dissection and basic model}
Slow/fast systems can be written in the form
\begin{equation}
\begin{split}
\dot{z} &= f(z,y,\e)\\
\dot{y} &= \e g(z,y,\e)
\end{split}
\label{fast-time}
\end{equation}
where $z$ is a vector of \textit{fast} variables, $y$ is a vector of \textit{slow} variables and $\e$ is the slow/fast timescale ratio.  Such systems arise in many areas of mathematical modeling and can describe general multiple timescale phenomena. In the singular limit, where $\e \rightarrow 0$, one obtains an equation where the slow variable(s) $y$ can be considered as parameter(s) for the fast subsystem ($z$). This approach allows one to investigate the dynamics of the fast subsystem and subsequently to construct solutions of the full system by carefully reintroducing slow dynamics -- i.e., by studying the perturbed system ($\e \neq 0$).  An elegant mathematical toolset known as geometric singular perturbation theory has been developed to study such phenomena~\cite{Cronin:1999p2714, Jones:2000p2718}.
 
The results brought forward in this paper relate to the effect of pulsatile perturbations on slow/fast bursters in which a delayed Hopf bifurcation is central to the onset of rapidly varying dynamics. One of the first systems of this type to be studied was a bursting phenomenon in the Belousov-Zhabotinskii chemical reaction~\cite{Rinzel:1982p10049}. The mechanisms central to our study can be found in many chemical, physical and biological systems. 

However, our main purpose is to better understand the effect of perturbations on conductance based models of single neurons (i.e., models of Hodgkin-Huxley type) that possess such a separation of timescales. In many cases, calcium concentration acts as a slow variable while voltage and associated ionic currents evolve on the fast timescale. Rinzel and Lee first studied fast/slow solutions to models of parabolic bursters in 1987 using singular perturbation methods \cite{Rinzel:1987p2710}. Since then, much effort has been invested in understanding bursting solutions arising in conductance-based neural models and their reduced forms~\cite{Terman:2005p2767,Rinzel:1989p2766}. In 2000, Izhikevich produced a classification of bursting mechanisms \cite{Izhikevich:2000p1909}, including all of the possible codimension one bifurcations of the fast subsystem that could be responsible for the onset and termination of spiking dynamics.
\begin{figure}[h]
\begin{center}
\includegraphics{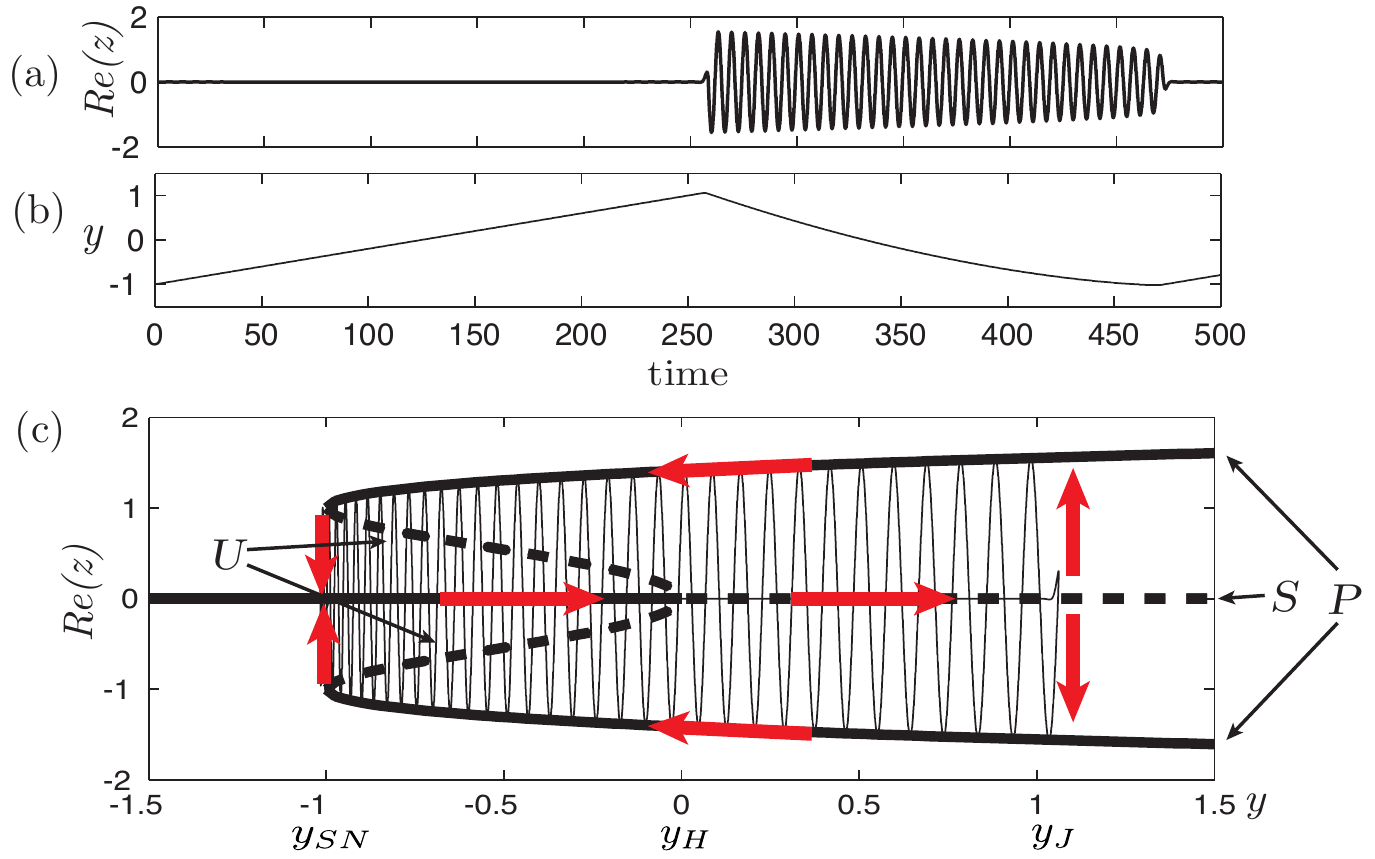}
\caption{Elliptic bursting trajectory from model (\ref{canon_bautin}). (a) Trace of $Re(z)$. (b) Trace of $y$. (c) Projected into the ($Re(z)$, $y$) plane, thicker black lines show $y$-parametrized fixed points and limit cycles of the fast subsystem (i.e., from the singular limit).  S is a family of fixed points (stable for the solid line and unstable for dashed), U is a family of unstable orbits and P is a family of stable orbits. Thin black line shows bursting trajectory, with red arrows indicating time evolution direction.}
\label{fig:bif}
\end{center}
\end{figure}

We concentrate on cells characterized as (subcritical) elliptic bursters or type III bursters, in which a subcritical Hopf bifurcation of the fast subsystem drives the onset of bursting and a saddle node bifurcation of limit cycles is responsible for the return to silent state (see Fig. \ref{fig:bif}). This type of intrinsic bursting is well understood in the absence of forcing (input).  In particular, Su, Rubin and Terman established existence and stability properties of elliptic bursting solutions in~\cite{Su:2004p2369}. Izhikevich \cite{Izhikevich:2000p2513} presented a portrait of elliptic bursting dynamics by describing the fast subsystem via the normal form of the (codimension two) Bautin bifurcation, and derived a closely-related model of weakly coupled networks of these bursters~\cite{Izhikevich:2001p2437}. Throughout this paper, we keep our analysis as general as possible and often illustrate our results with a variation of Izhikevich's \textit{normal form model} (see Eqn.\eqref{canon_bautin}) to simplify mathematical manipulations. However, we stress that our analysis can be carried out for any elliptically bursting model. 

The \textit{normal from model} is
\begin{equation}
\begin{split}
\dot{z}&=(y+iw)z+2z|z|^2-z|z|^4+I(t)\\
\dot{y}&=\e(a-|z|^2-by)
\end{split}
\label{canon_bautin}
\end{equation}
where $z \in \C$ represents the fast variable and $y \in \R$ the slow variable. This is as found in~\cite{Izhikevich:2001p2437}, but with the term $by$ added to the slow variable dynamics in order to explore the effect or nonlinear ramping in the delayed Hopf bifurcation, a feature that is found in many models for which bursting dynamics are caused by slowly varying calcium concentration.  We set $w=1$ and $\e=0.01$ for the remainder of the paper and will consider distinct cases in which we vary $a$ and $b$. Intrinsically bursting solutions generally arise for parameter choices where $\dot{y}$ is positive when $z$ is at rest and negative when $z$ is spiking (oscillating).  

Notice the forcing via the signal $I(t)$ in the equation for the fast variable. We wish to model an input signal that causes an instantaneous voltage response, so we set $I(t)=\sum_n A_n \delta(t-t_n)$ where $\delta$ is a delta function and $A_n \in \R^+$. We will refer to these perturbations as $t_n$-\textit{kicks} of \textit{amplitude} $A_n$: a kick simply translates a solution at time $t_n$ by an amount $A_n$ in the real ``direction" of $z$. In this paper, we focus on periodic kicks of fixed amplitude $A=A_n$ and equal spacing $\tau=t_n-t_{n-1}$ between kick times.

\subsection{The elliptic bursting cycle}

Figure \ref{fig:bif} shows an elliptic bursting trajectory from numerically integrating Eqn.~(\ref{canon_bautin}) with $a=0.8$ and $b=0$.  As mentioned above, the standard approach is to think of the slow variable $y$ as a parameter that determines the fast dynamics ($z$). The fast subsystem undergoes a subcritical Hopf bifurcation at $y_H=0$ and a saddle-node of limit cycle bifurcations at $y_{SN}=-1$. When $\e=0$, and given a fixed $y \in (y_{SN},y_H)$, the fast subsystem is bistable and has a sink at $z=0$ inside an unstable periodic orbit, itself contained within a stable periodic orbit. These $y$-parametrized limit sets form normally hyperbolic invariant manifolds that persist under small perturbations ($1>>\e >0$)~\cite{Su:2004p2369}. 
We denote by $S$ the family of equilibria ($z=0$).  These are stable when $y$ lies to the left of $y_H$ and unstable on the right, together forming the \textit{silent} branch.  $P$ is the family of stable periodic orbits and represents the \textit{spiking} branch. For Eqn.\eqref{canon_bautin}, the radii of these orbits about $S$ are given by $r_P(y)=\sqrt{1+\sqrt{y+1}}$ for $y > y_{SN}$. Finally, $U$ refers to the family of unstable orbits with radii $r_U(y)=\sqrt{1-\sqrt{y+1}}$ for $y \in (y_{SN},y_H)$, acting as a separatrix between the stable side of $S$ and $P$. A bursting solution occurs when the slow dynamics of $y$ steer the fast subsystem rightward along the silent branch $S$, until the Hopf point $y_H$ is reached. The solution then keeps moving rightward, sticking close to $S$ for some transient period even though the equilibria forming $S$ are now unstable, but is eventually attracted to $P$ when $y$ reaches $y_J$, where spiking begins. The $y$ dynamics then pull the oscillating fast subsystem leftward along $P$, until the latter vanishes at $y_{SN}$, where the solution is attracted back to $S$ and another cycle begins. 

Of particular interest is the {\it slow passage effect} through the Hopf point $y_H$, in which the solution does not immediately jump up into spiking when $S$ looses stability ($y_J \neq y_H$). This delayed bifurcation phenomenon as been previously studied~\cite{Baer:2008p11521,Baer:1989p4834,Su:2004p2369} and its implications in the context of pulsatile perturbations will be established in Sect. \ref{map}. Importantly, several authors have shown that noise can sharply diminish this slow passage effect \cite{Kuske:2002p11725,Kuske:1999p5231,Baer:1989p4834,Su:2004p2369}. While we first treat the noiseless case, we study the effect of stochastic terms on the phase reduction and response dynamics in Sect.~\ref{noise}. 


\subsection{Phase reduction}
\label{redux}

There is a substantial body of literature concerning phase reduction of oscillators and their behavior under noise, forcing, or coupling~\cite{EK84,Izhikevich:2007p1301,Brown:2004p122,Croisier:2009p2491,Golubitsky:2006p150,Nakao:2005p126,Terman:2008p153,Ermentrout:2010p10442}. The general idea is simply to associate the endpoints of a periodic solution's cycle and to parametrize the movement along this solution with a phase $\o \in S^1$. Although the point $\o=0$ is arbitrary, it is often chosen to correspond to a distinguishable event within the periodic orbit, such as the apex of an action potential in a model of a spiking cell. This reduction becomes useful when the limit cycle has some stability properties and one can track the phase response of the solution following a perturbation:  specifically, by computing the phase difference on $S^1$ between the unperturbed solution and the perturbed one, as $t \to \infty$ and the latter contracts back to the limit cycle. When well defined, this one-dimensional description has the advantage of being analytically tractable while preserving the behavior of an oscillator subject to perturbations.

For systems with asymptotically stable limit cycles, phase reduction can be carried out rigorously as long as the effect of forcing translates the solution to a point in the cycle's basin of attraction. This basin is foliated by the strong stable manifolds of each point on the limit cycle. By knowing on which manifold -- or {\it isochron} -- the solution lands following a perturbation, we know exactly to which phase it will be attracted in the limit as $t \to \infty$ \cite{EK84,Winfree:2001p9818,G75,Medvedev:2001p4516,Croisier:2009p2491,Izhikevich:2000p2936}.

We next describe a phase reduction for elliptic bursters.  As we will discuss further, the bursting trajectories are not stable limit cycles, and so we cannot directly compute isochrons. Nevertheless, the difference between the timescale of the burst period and the timescale of attraction normal to the (singular limit) solution enables us to proceed. We closely follow Best et al. \cite{Best:2007p128}, who derive circle dynamics for an elliptic burster receiving periodic inputs from a model excitatory neuron. In their work, each excitatory kick always transitions a solution from the resting to the spiking state -- i.e., to the branch of periodic orbits $P$ in Fig.~\ref{fig:bif}, if it is not already following that branch.  Best et al. also use an approximation of linearity for the slow trajectories (i.e., $\dot y$ is piecewise constant).   Here, we relax both of these assumptions, in particular while studying the response to weaker kicks that do not necessarily generate a burst.  As we will see, it is these weaker kicks that will lead to desynchrony; that is, dynamics that appear chaotic or are phase locked at high period.

In \cite{Su:2004p2369}, the authors use Fenichel theory to show that there exist $O(\e)$ neighborhoods  $N_S$ and $N_P$ around $S$ and $P$ such that -- if a solution enters either the left side of $N_S$ or the right side of $N_P$ and the slow dynamics behave as mentioned above -- then the solution will transition between the two neighborhoods in a periodic, bursting fashion. Furthermore, they use averaging techniques to show that the period of such a cycle can be approximated up to $O(\e)$ by the sum of the passage times $T_S$ and $T_P$ along the respective manifolds.

Although it is not clear whether or not there exists a single periodic solution of the full system, \cite{Su:2004p2369} shows that such dynamics are at least metastable. That is, $N_S$ and $N_P$ are locally attracting and any bursting solutions must live in these sets. Furthermore, solutions starting outside the two attracting sets are quickly attracted back to them.
 
Numerically obtained solutions of Eqn. (\ref{canon_bautin}) do not trace back exactly the same path from cycle to cycle but have periods that vary minimally, as expected~\cite{Su:2004p2369}.  Specifically, we numerically integrated a solution in order to obtain 150 burst cycles and computed the coefficient of variation ($CV=\frac{\textrm{standard deviation}}{\textrm{mean}}$) of the cycle durations.  We found $CV=O(10^{-3})$ for parameter choices yielding bursting solution in Eqn. \eqref{canon_bautin} (numerical methods as in Sect.~\ref{numerics}).  Whether there is a periodic solution with a much longer period than the bursting cycle, or whether solutions are instead quasiperiodic or aperiodic remains an open question. 
 
The small $CV$ (indicating a robust cycle period), along with the metastability described above, motivate an approximate reduction to dynamics on the circle, as in~\cite{Best:2007p128}. We will revisit the notion of uncertainty in cycle periods in Sect. \ref{noise}. However, for what follows, we use the singular limit assumption that a bursting trajectory evolves along $S$ and $P$ with well defined passage times $T_S$ and $T_P$, and will use this trajectory to compute phase reduced dynamics.

\begin{figure}[h]
\begin{center}
\includegraphics[width=350pt]{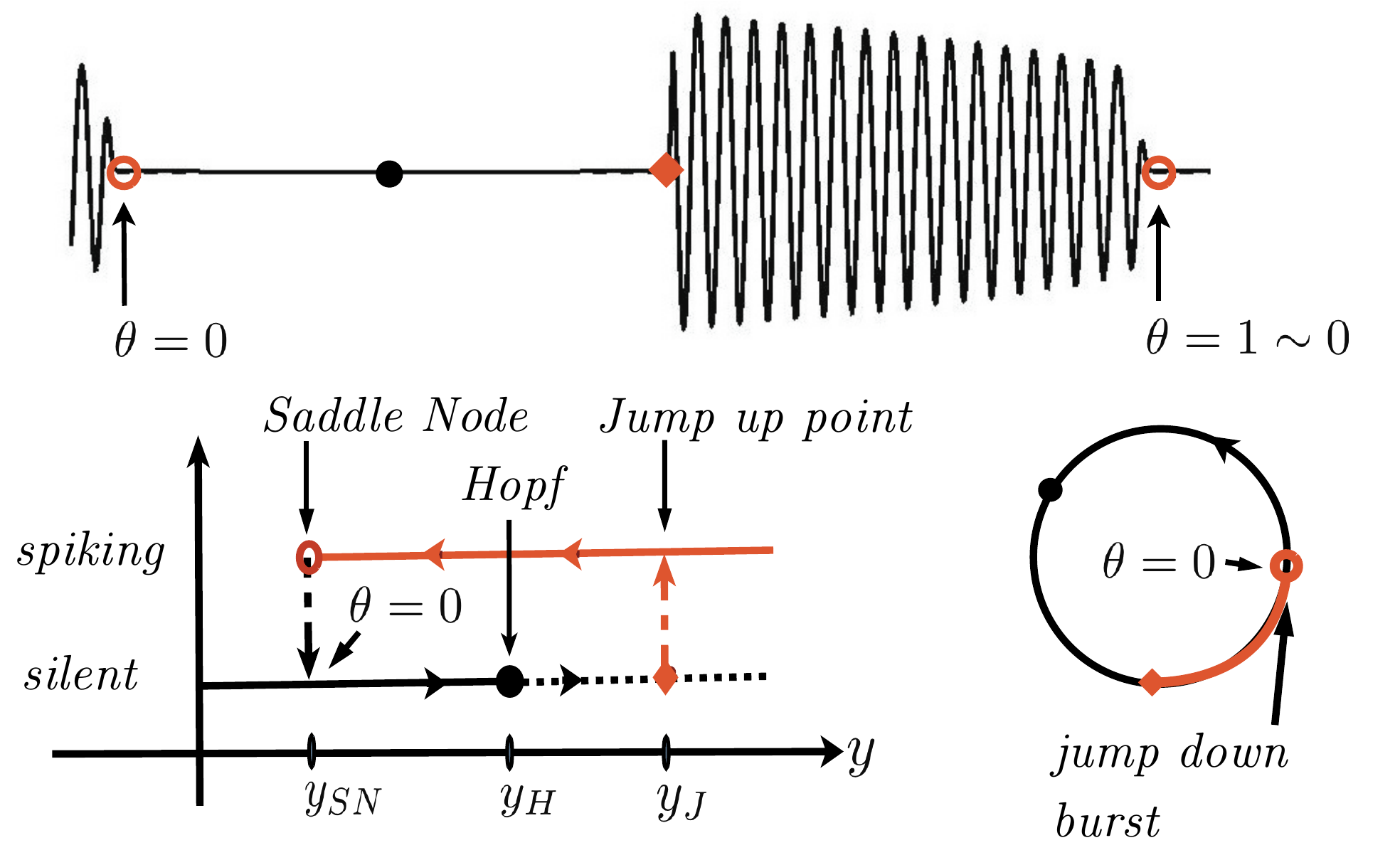}
\caption{Schematic presentation of phase reduction. We associate endpoints of a bursting cycle (where $y=y_{SN}$), and map the trajectory onto the unit circle.}
\label{fig:redux}
\end{center}
\end{figure}
 
As illustrated in Fig. \ref{fig:redux}, we represent bursting solutions by a phase variable $\o(t) \in S^1=\R/\Z$.  We let $\o=0\simeq1$ correspond to the ``jump down point" on the bursting trajectory, where solutions transfer from spiking to resting. We choose this reference point because the spiking to resting transition is fast and is associated with a constant value of the slow variable, $y=y_{SN}$ (unlike, as we will see, the transition to spiking following a pulsatile input). Essentially, the phase $\o$ of a bursting cycle is given by time rescaled by the period ($\o=\frac{t}{T_S+T_P}$) where at $t=0$, $y=y_{SN}$. Note that the first portion of the unit circle following $\o=0$ represents the silent branch $S$ and the remaining portion represents the spiking branch $P$. 

\subsection{The kick map}
\label{map}

We now study phase dynamics of bursters receiving pulsatile inputs (kicks).   
Specifically, we derive a phase translation mapping $F_A(\o)$, such that if $\o$ is the phase of a cell when a kick of strength $A$ arrives, $F_A(\o)$ is the phase of the kicked solution -- relative to the unperturbed solution -- after it relaxes back to the burst cycle.  We will refer to this as the \textit{kick map}.  In \cite{Best:2007p128}, a similar map for elliptic bursters is derived, and this idea inspired the present work.  However, there are two differences with the map we derive here.  First, the phase of trajectories in~\cite{Best:2007p128} is defined relative to the period of pulsatile inputs; in our case, phase is defined relative to the (unperturbed) period of a burst trajectory.  This latter construction has been previously used in the context of integrate and fire cells with soft reset \cite{Belair:1986p8707}, and we find that this makes it easier to visualize the role of changing kick amplitude and period on the structure of the map.  Second, in~\cite{Best:2007p128}, only ``strong" kicks are considered; as we will see, the map develops additional features, including discontinuity and expansion, in the case of weaker kicks.

As discussed in Sect.~\ref{redux}, we assume that unperturbed elliptic bursting solutions have fixed times spent in silent ($T_S$) and spiking phases ($T_P$).  When computing a map for a given system as done in Sect. \ref{canon_map}, it is best to work with unscaled time and later (implicitly) normalize the time variable by the burst period ($T_S+T_P$) so that our phase variable $\o$ remains between zero and one. In this section however, we derive the kick map for an arbitrary elliptic bursting model and assume that the period is already unitary ($T_S+T_P=1$) for simplicity. The rest of the notation follows that of system~\eqref{canon_bautin} but the reader should keep in mind that $z$ and $y$ represent general fast and slow variables.

Our computations are intimately linked to the evolution of the slow variable $y$ along the branches $S$ and $P$. Recall that $y$ spans $[y_{SN},y_J]$ when $z$ travels along $S$ or $P$. To better track the variable $y$, we label its dynamics along $S$ by $y(t)=h_S(t)$ and along $P$ by $y(t)=h_P(t)$. 
 
Thus, for a burst trajectory that starts at $y=y_{SN}$ when $t=0$, we have\begin{equation}
y(t) = \left\{ \begin{array}{ll}
h_S(t) & \textrm{if $0\leq t < T_S$}\\ 
h_P(t) & \textrm{if $T_S \leq t < 1$}
\end{array} \right.
\label{y_dyn}
\end{equation}
where $h_S(0)=y_{SN}=h_P(1)$ and $h_S(T_S)=y_J=h_P(T_S)$. Here, $h_S$ and $h_P$ are functions with ranges $[y_{SN},y_J]$ and respective domains $[0,T_S]$ and $[T_S,1]$. We assume $h_S$ is an increasing function while $h_P$ is decreasing. We now define the phase $\o$ of a bursting solution $(z(t),y(t))$ by
\begin{equation}
\o = \left\{ \begin{array}{ll}
h_S^{-1}(y(t)) & \textrm{if silent ($z(t) \in S$)}\\ 
h_P^{-1}(y(t)) & \textrm{if spiking ($z(t) \in P$)}.
\end{array} \right.
\label{phase_def}
\end{equation}
For unperturbed solutions, $\o=t \mod 1$. In general, expressions for $h_S$ and $h_P$ can be hard to find.  As in standard approaches, one can integrate the dynamics of $y$ by restricting the fast variable to the invariant manifolds $S$ and $P$, and using the averaged motion of $z$ on those manifolds~\cite{Rinzel:1989p2766, Ermentrout:2010p10442,Baer:1989p4834,Su:2004p2369}.  These calculations yield explicit formulas for the normal form system (\ref{canon_bautin}), as we demonstrate in Sect. \ref{canon_map}.

\begin{figure}[h]
\begin{center}
\includegraphics{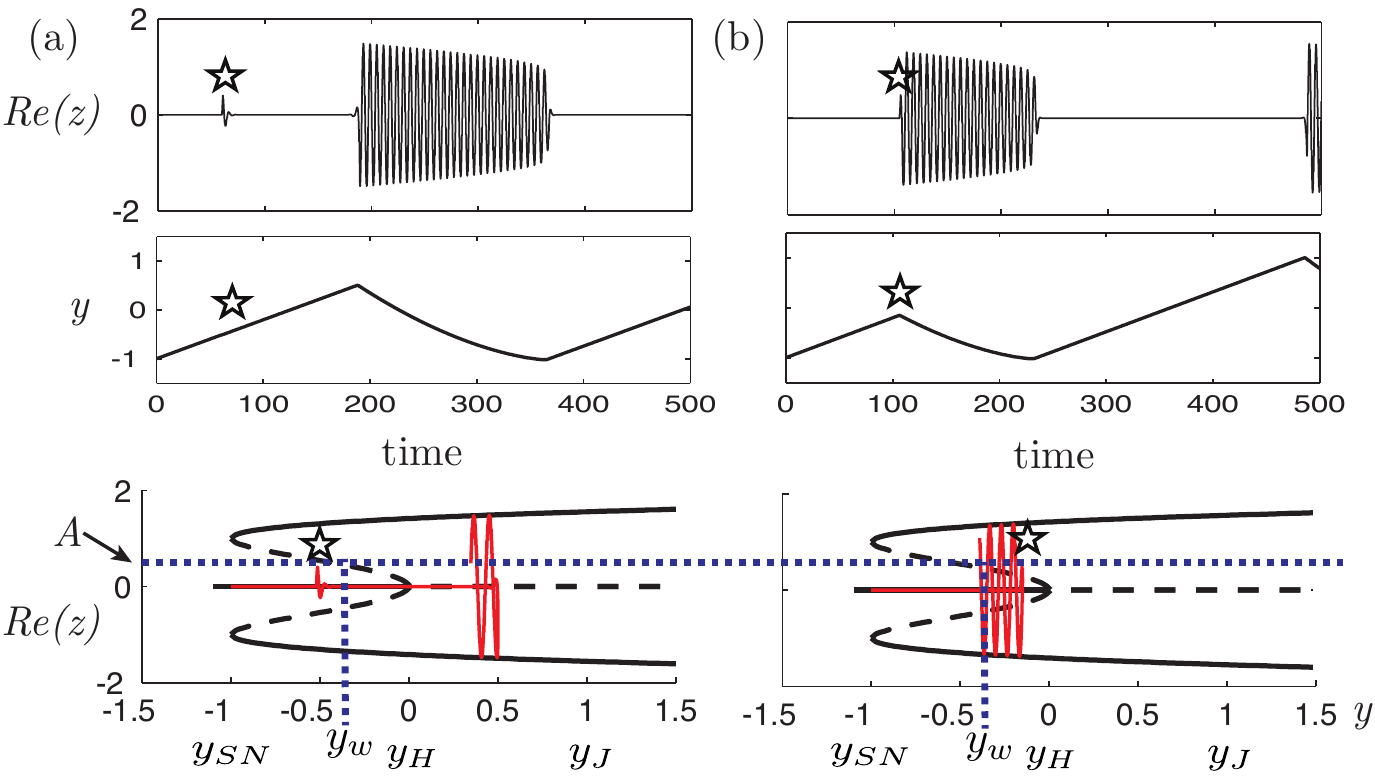}
\caption{Weakly kicked ($A=0.5$) trajectories from (\ref{canon_bautin}). Top to bottom:  $Re(z)$ trace, $y$ trace and $y-Re(z)$ solution curve (red) truncated at $t=200$ for clarity. Star indicates a kick and dashed blue line the kick's amplitude $A$. (a) Kick received at $t=80$ and does not clear separatrix. (b) Kick received at $t=110$ and clears separatrix.}
\label{fig:geometry}
\end{center}
\end{figure}

We are now equipped to define the kick map.
Our first assumption is that a kick that arrives during the spiking regime has no effect on $\o(t)$. Due to the separation of timescales, trajectories are attracted back to the stable limit cycles that form $P$ in a vanishingly short time (with respect to the timescale over which the phase evolves).  On the other hand, a kick while the cell is silent ($z$ near $S$) may have distinct outcomes. If $y(t) \in  [y_H,y_J]$ (sticking to the unstable part of $S$), any kick will send the cell to the spiking state since the trajectory is highly sensitive to perturbations. However, if $y(t) \in [y_{SN},y_H]$, one of two things can happen. If the kick is strong enough to translate the solution past the separatrix $U$, then the cell jumps to the spiking state (Fig. \ref{fig:geometry} (b)). If on the other hand, the kick is not strong enough, the solution will attract back to $S$ (Fig. \ref{fig:geometry} (a)). 
\beq
\o_w(A)=h_S^{-1}(y_w(A)).
\label{cutoff}
\eeq

From now on, we refer to a kick as \textit{strong} if $A>r_U(y)$ for all $y \in [y_{SN},y_H]$ where, $r_U(y)$ is the distance between $S$ and $U$ at $y$, in the direction of the kick.
 In other words, a strong kick will immediately result in spiking independently of the cell's phase. In contrast, we define a \textit{weak} kick as one with an amplitude $A<r_U(y)$ for values of $y$ in some subinterval of $[y_{SN},y_H]$, so that the kick does not always immediately cause a cell to spike. (Note that a weak kick will result in spiking for any value of $y$ outside this interval.)

\begin{figure}[h]
\begin{center}
\includegraphics[width=350pt]{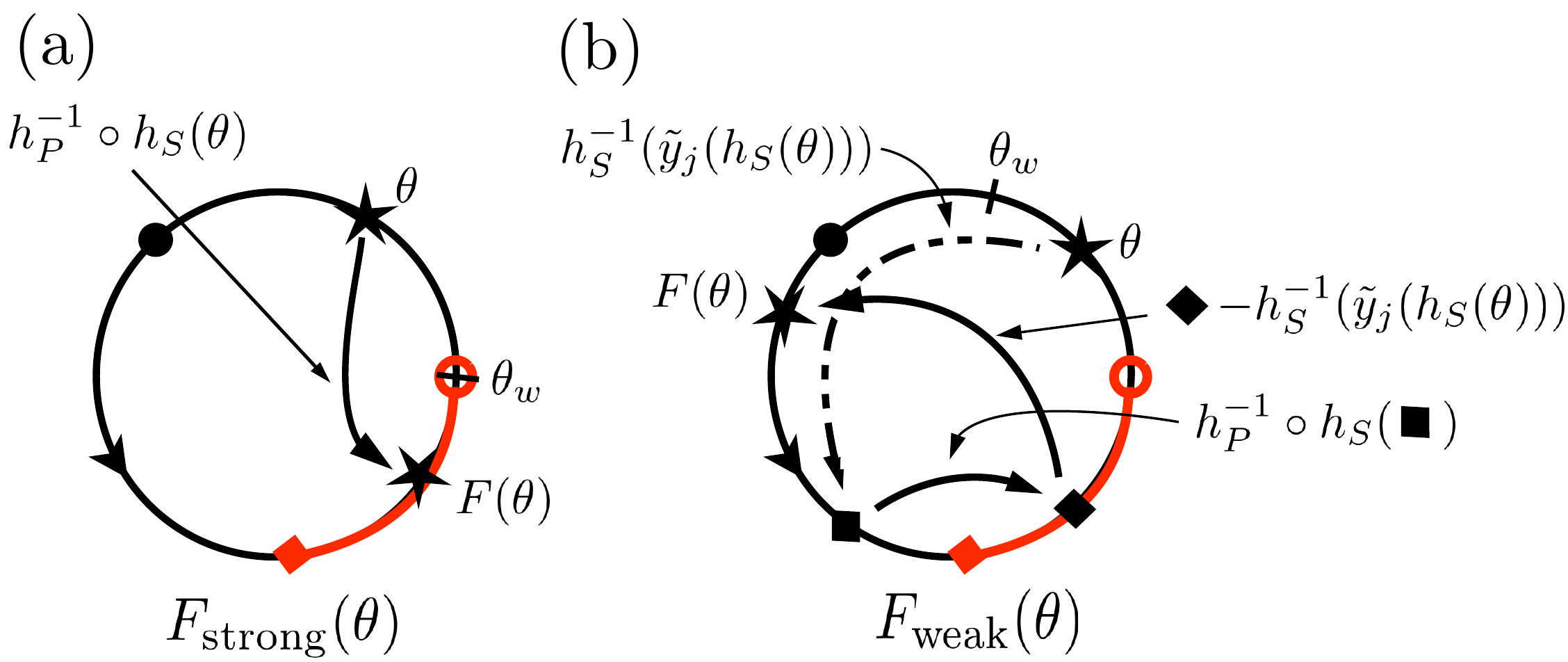}
\caption{Schematic representation of the kick map $\o \to F_A(\o)$ on the unit circle. (a) Strong kick map (\ref{strong_map}). (b) Weak kick map of Eq.~(\ref{weak_map}). The solid arrows represent instantaneous change of phase while the dotted arrow indicates the evolution of the phase along the circle in time.} 
\label{fig:diag_map}
\end{center}
\end{figure}

For the strong kick case, the kick map is\begin{equation}
F_{A}(\o) = \left\{ \begin{array}{ll}
h_P^{-1}\circ h_S(\o) & \textrm{if $\o \in [0,T_S]$}\\ 
\o & \textrm{if $\o \in [T_S,1]$}.
\end{array} \right.
\label{strong_map}
\end{equation}
While in silent phase, $\o \in [0,T_S]$, $y$ increases according to $h_S(t)$; when spiking is induced by a kick, the value of $y$ is left unchanged but its dynamics are ``reversed" and it starts decreasing via $h_P(t)$.  The cell will then spike until $y$ reaches $y_{SN}$, which takes less time since we started closer $y_{SN}$.  Thus the impact of the kick is to advance the phase.  This explains the first line of Eqn.~\eqref{strong_map}, and is sketched in the panel (a) of Fig. \ref{fig:diag_map}.  As already discussed, the kick has no effect when the cell is spiking, as expressed in the second line of  Eqn.~\eqref{strong_map}.

For weak kicks, the situation is more complex. Recall that the branches $U$ and $P$ meet at $y_{SN}$ ($r_U(y_{SN})=r_P(y_{SN})$), and  $U$ vanishes at $y_H$ ($r_U(y_H)=0$); we assume that $r_U$ is a non-increasing, continuous function. As a result, for $r_U(y_{SN})>A>0$, we can find $y_w(A)$ such that $A=r_U(y_w(A))$ and $A<r_U(y)$ for all $y \in [y_{SN},y_w(A)]$.
If the cell is in silent phase, $y_w(A)$ is essentially a cutoff point before which a weak kick will not elicit immediate spiking, as illustrated in Fig. \ref{fig:geometry}(a).  A weak kick delivered at any other point through the burst cycle will result in instantaneous spiking, as in the strong kick case. We recast this condition in phase coordinates, obtaining the cutoff phase

What happens to a cell's phase when a weak kick does not evoke spiking ($\o \in [0,\o_w(A)]$)?  The trajectory is attracted back toward $S$ but jumps up into spiking before it reaches $y_J$, as if it retained a memory of this past weak kick (Fig. \ref{fig:geometry} (a)). To better understand this phenomenon, we first need expressions for slow passage times through Hopf points. We use results for delayed bifurcations derived in~\cite{Baer:2008p11521} in order to predict points of transition between silent and active states.

Let $\l(y)$ be the extremal eigenvalue of the fast subsystem linearized about the equilibrium points $z_0(y)$, for some chosen value of $y \in [y_{SN},y_J]$.  (Recall that this collection of points forms the manifold $S=\{z_0(y) | y \in [y_{SN}, y_J]\}$.)  Points on $S$ to the left of $y_{H}$ are sinks with $Re(\l(y))<0$, whereas to the right, they are sources and $Re(\l(y))>0$. Assume $\frac{d\l(y)}{dy}|_{y=y_H} \neq 0$. As a solution is pulled to the right by the slow dynamics, $y$ crosses $y_H$, $Re(\l(y))$ changes sign and solutions switch from being attracted to being repelled by $S$. However, this repulsion is not immediately apparent:  the difference in time scales of our slow/fast system induces a discrepancy in spatial scales. As $y$ varies, $z$ is attracted to $S$ for $y<y_H$ and repelled for $y>y_H$, both at an exponential rates (proportional to $Re(\l(y))$) on the fast timescale. 

The length of the slow passage (also called delay) to the right of $y_H$ depends on the dynamics of $y$ along $S$. Borrowing notation from~\cite{Baer:2008p11521}, suppose we can write
\beq
y(t)=y_i+g(\e t)
\label{slow}
\eeq
where we assume $g$ is a non-decreasing function and $g(0)=0$. In this context, $h_S(t)=y_{SN}+g(\e t)$. If the system's full solution starts at $(z_i,y_i)$ such that $y_i \in [y_{SN},y_H]$ and $z_i$ is far enough from  $z_0(y_i)$ but is still in its basin of attraction, then the jump up point $y_j$ can be implicitly computed via
\begin{equation}
0=\int_{y_i}^{y_j} (\frac{d}{dy}g^{-1}(y-y_i))Re(\l(y))dy
\label{integral_condition}
\end{equation}
where $g^{-1}$ is the inverse of $g$. In the unperturbed case, the point $y_J$ can be derived using $y_i=y_{SN}$ in (\ref{integral_condition}); this is sometimes called the memory effect for elliptic bursters~\cite{Rinzel:1987p2710}. We refer the reader to the appendix of~\cite{Baer:2008p11521}  for the derivation of this integral condition.

If a weak kick does not elicit instantaneous spiking, we find that the ``memory" starts anew at the time of the kick; in other words, if the kick is administered when $y\in [y_{SN},y_w(A)]$,  then we set $y_i=y$ in Eqn.~(\ref{integral_condition}), and denote the associated jump up point by $\tilde{y}_j(y)$.  In phase coordinates, if the cell is kicked at $\o \in [0,\o_w(A)]$, $y$ is given by $h_S(\o)$. The onset of spiking happens at $y_j=\tilde{y}_j(h_S(\o))$, or equivalently at $\o_j=h_S^{-1}(\tilde{y}_j(h_S(\o)))$. 
We capture this via:
\begin{equation}
F_{A}(\o) = \left\{ \begin{array}{ll}
\o+h_P^{-1}(\tilde{y}_j(h_S(\o)))-h_S^{-1}(\tilde{y}_j(h_S(\o))) & \textrm{if $\o \in [0,\o_w(A))$}\\
h_P^{-1}\circ h_S(\o) & \textrm{if $\o \in [\o_w(A),T_S]$}\\ 
\o & \textrm{if $\o \in [T_S,1]$}.
\end{array} \right.
\label{weak_map}
\end{equation}
The first conditional definition maps the slow variable to its jump up value $y_j=\tilde{y}_j(h_S(\o))$, and then applies the strong kick map $h_P^{-1}\circ h_S(\o)$, and finally translates back by the phase $-h_S^{-1}(\tilde{y}_j(h_S(\o)))$ to account for the time taken in ``slow passage" from $y_i=h_S(\o)$ to $y_j=\tilde{y}_j(h_S(\o))$. This mapping is sketched in the right panel of Fig. \ref{fig:diag_map}; 
note that it is only valid if the cell does not receive additional kicks before it enters the spiking state.

This construction implies that the shape of a kick map does not vary continuously with the strength of a kick. Moreover, the only criterion which dictates the qualitative shape of the map is the value of $\o_w(A)$ (we sometimes drop the $A$ and write $\o_w$).
Other than determining the threshold $\o_w$, the perturbative role of $A$ is dimensionless since a kick acts on the fast variable as opposed to the phase which is defined over the slow timescale. That is, if two kicks of distinct amplitudes yield the same strong (resp. weak) outcome, the discrepancies between the times needed to attract close to the unperturbed trajectories are negligible. We emphasize that $\o_w$ decreases as $A$ increases:  for strong kicks with $A>r_U(y_{SN})$, $\o_w$ is always zero and the map does not change shape as $A$ increases further.
We use Eqn. (\ref{weak_map}) as the general expression for our kick map. 

We stress that there is a fundamental difference between the maps induced by strong and weak kicks.  A weak kick always induces expansion in the kick map, as long as it speeds passage ($\tilde{y}_j \neq y_J$) through the Hopf point. Indeed, notice that by construction, $h_P^{-1}(y) > h_S^{-1}(y)$ for any $y \in (y_{SN},y_J)$ and $h_P^{-1}(y) = h_S^{-1}(y)$ when $y=y_{SN},y_J$. From expression (\ref{weak_map}), we see that $F_A(\o)>\o$ on $(0,\o_w)$ and $F_A(0)=0$. It follows from the mean value theorem that $\frac{dF_A}{d\o}>1$ on some region contained in $[0,\o_w)$. We also note that when $h_P$ and $h_S$ have similar shapes, it generally implies expansion of $F_A$ on the whole interval $[0,\o_w)$.
To better illustrate this and other features, we compute expressions of this map for system (\ref{canon_bautin}).

\subsection{Computing the kick map for the normal form model}
\label{canon_map}
In this section, we derive an analytical approximation of the kick map for the elliptic bursting normal form model. The task at hand is simple:  use Eqn. (\ref{canon_bautin})  to compute the ingredients of expression (\ref{weak_map}):  $h_S(\o)$, $h_S^{(-1)}(\o)$, $h_P^{-1}(\o)$, $\tilde{y}_j(y)$ and $\o_w(A)$. 

We first turn to $h_S$ and $h_P$, which are essentially the $y$ components of a solution to Eqn. (\ref{canon_bautin}), in silent and spiking modes respectively. In contrast with the previous section, we carry out computations using unscaled time which implies a full burst period $T \neq 1$. One can still think of $t$ as $\o$ in what follows, but the expression of the final map has to be rescaled.

We exploit the separation of timescales in our equation and make the assumption, as in the singular limit, that 
$z$ evolves exactly on the manifolds $S$ and $P$. Notice that the $y$ dynamics depend linearly on $|z|^2$. By substituting $|z|^2=0$ for $h_S$, and $|z|^2=r_P(y)$ for $h_P$, we obtain two scalar O.D.E.s
\begin{eqnarray}
\dot{y}= \e( a -by)& \quad \textrm{when $z$ is on $S$}\label{S_eq}\\
\dot{y}= \e (a-r_P(y)^2-by) & \quad \textrm{when $z$ is on $P$}. 
\label{P_eq}
\label{ODE_uncoupled}
\end{eqnarray}
Bursting occurs when the right hand sides of \eqref{S_eq} and \eqref{P_eq} are respectively positive and negative, steering the fast dynamics in the required directions along $S$ and $P$. Here we concentrate on parameter values $a>0$ and $b \geq 0$ which satisfy this condition.
As we will shortly see, $b>0$ implies that the evolution of $y$ along $S$ follows a saturating exponential ramp while $b=0$ implies a linear ramp. Both scenarios are found in biological systems, and -- as we now show -- the resulting kick maps have common characteristics which unify their response to pulsatile perturbations.  

For Eqn.  \eqref{S_eq}, we can easily solve and get
\begin{eqnarray}
h_S(t)=\e at+C_S & \quad b=0 \label{hs0}\\
h_S(t)=\frac{ a}{b}+C_Se^{-\e bt} & \quad b>0 \label{hs}
\end{eqnarray}
where $C_S$ is an integrating constant. Setting $y=y_{SN}$(=$-1$) at $t=0$, it is easy to see that $C_S=y_{SN}$ for $d=0$ and $C_S=y_{SN}-\frac{ a}{b}$ for $b>0$. In turn, we have
\begin{eqnarray}
h_S^{-1}(y)=\frac{y-y_{SN}}{\e a} & \quad b=0 \label{inv_hs0}\\
h_S^{-1}(y)=-\frac{1}{\e b} \ln(\e[a-by])+\frac{1}{\e b}\ln(\e[a-by_{SN}]) & \quad b>0. \label{inv_hs}
\end{eqnarray}

For Eqn. \eqref{P_eq}, recall that 
$r_P(y)=\sqrt{1+\sqrt{y+1}}$ when $y>y_{SN}$. Solving this O.D.E. is not as straightforward; fortunately, we only need the inverse $h_P^{-1}$ (which is a proxy for $t$ in this context) to compute our map.  This can be obtained directly via integration:
\begin{eqnarray*}
h_P^{-1}(y)=-\frac{2}{\e}[(a-1)\ln(-a+\sqrt{y+1}+1)+\sqrt{y+1}]+C_P & \,b=0 \label{inv_hp0}\\
h_P^{-1}(y)= \frac{1}{\e b}\left[\frac{2\tan^{-1}\left(\frac{2b\sqrt{y+1}+1}{\sqrt{-4(a-1)b-4b^2-1}}\right)}{\sqrt{-4(a-1)b-4b^2-1}} -\ln(-a+by+\sqrt{y+1}+1)\right]+C_P& \, b > 0 \label{inv_hp}.
\end{eqnarray*}
 Using \eqref{inv_hs0}, \eqref{inv_hs} it is straightforward to compute the time ($T_S$) it takes for $h_S(t)$ to reach $y_J$, and then derive values for $C_P$ such that $h_P^{-1}(y_J)=T_S$ as required by our definition of $h_p$.

 To compute the jump up point $y_J$, we derive an expression for $\tilde{y}_j(y_i)$ which is also needed in the definition of our kick map. In order to use the integral condition (\ref{integral_condition}), we must first write expressions for the $y$ dynamics in the form of Eqn. \eqref{slow}. Here we substitute $y_{SN}$ by $y_i$ in the expression of $C_S$ for \eqref{hs0} and \eqref{hs} to allow for arbitrary initial conditions and get
 \begin{eqnarray*}
 y(t)=y_i+a\e t & \quad b=0\\
 y(t)=y_i+(1-e^{-b\e t})(\frac{a}{b}-y_{i}) & \quad b>0
 \end{eqnarray*}
 which yields
 \begin{eqnarray}
g(\e t)=a \e t & \quad b=0 \label{fret0}\\
g(\e t)=(1-e^{-b\e t})(\frac{a}{b}-y_i) & \quad b>0 \label{fret}
\end{eqnarray}
which in turn give
\begin{eqnarray}
g^{-1}(y)=\frac{y}{a} & \quad b=0 \label{zut0}\\
g^{-1}(y)=-\frac{\ln(1+\frac{y}{y_i-a/b})}{b} & \quad b>0 \label{zut}.
\end{eqnarray}

It is easy to compute the extremal eigenvalue for the linearization about $z=0$ of the fast subsystem in Eqn. \eqref{canon_bautin} :  $\l=y\pm iw$ and therefore $Re(\l(y))=y$. Turning now to the integral condition \eqref{integral_condition}, when $b=0$ we use \eqref{zut0} and write
$$0=\int_{y_i}^{y_j}\frac{y}{a}dy$$
by which we can deduce that
\begin{eqnarray}
y_j=\tilde{y}_j(y_i)=-y_i & \quad b=0.
\label{yj0}
\end{eqnarray}
In other words, when the slow variable follows a linear ramp along the silent branch $S$ (see \eqref{fret0}),
 the jump up point $y_j$ is symmetric to the initial point $y_i$ about $y_H=0$. This can be seen from Fig. \ref{fig:bif} where $y_{SN}=-1$ and $y_J=1$.

In the case where $b>0$, $y$ follows a saturating exponential ramp along $S$ (see \eqref{fret}) which implies that $y$ decelerates as it moves rightwards and thus shortens the slow passage.  Using \eqref{zut} in the integral condition \eqref{integral_condition}, we get
$$0=\int_{y_i}^{y_j}\frac{y}{b(1+\frac{y-y_i}{y_i-a/b})(a/b-y_i)}dy$$
which gives
$$0=\left[-\frac{1}{b} \frac{a \ln(by-a)}{b}+y\right]_{y=y_i}^{y_j}$$
thus implying the relation
$$y_i-y_j=\frac{a}{b} \ln(\frac{by_j-a}{by_i-a}).$$
We then isolate $y_j$ to get
\begin{eqnarray}
y_j=\tilde{y}_j(y_i)=\frac{a}{b}[W(-\frac{1}{a} e^{\frac{b}{a}y_i-1}(a-by_i))+1] & \quad b>0
\label{yj}
\end{eqnarray}
where $W$ is the Lambert W function (product logarithm function).

The last ingredient we need is an expression for the cutoff value $y_w(A)$; recall that this marks the $y$ boundary below which a given kick will not clear the separatrix $U$ with radius $r_U(y)=\sqrt{1-\sqrt{y+1}}$. This quantity does not depend on the slow dynamics and is therefore the same for both our cases. For a kick amplitude $A$, $y_w(A)=(1-A^2)^2-1$ if $A \in [0,1]$ (a weak kick). When $A>1$, the kick is strong and we set $y_w=y_{SN}=-1$. 

\begin{figure}[t]
\begin{center}
\includegraphics{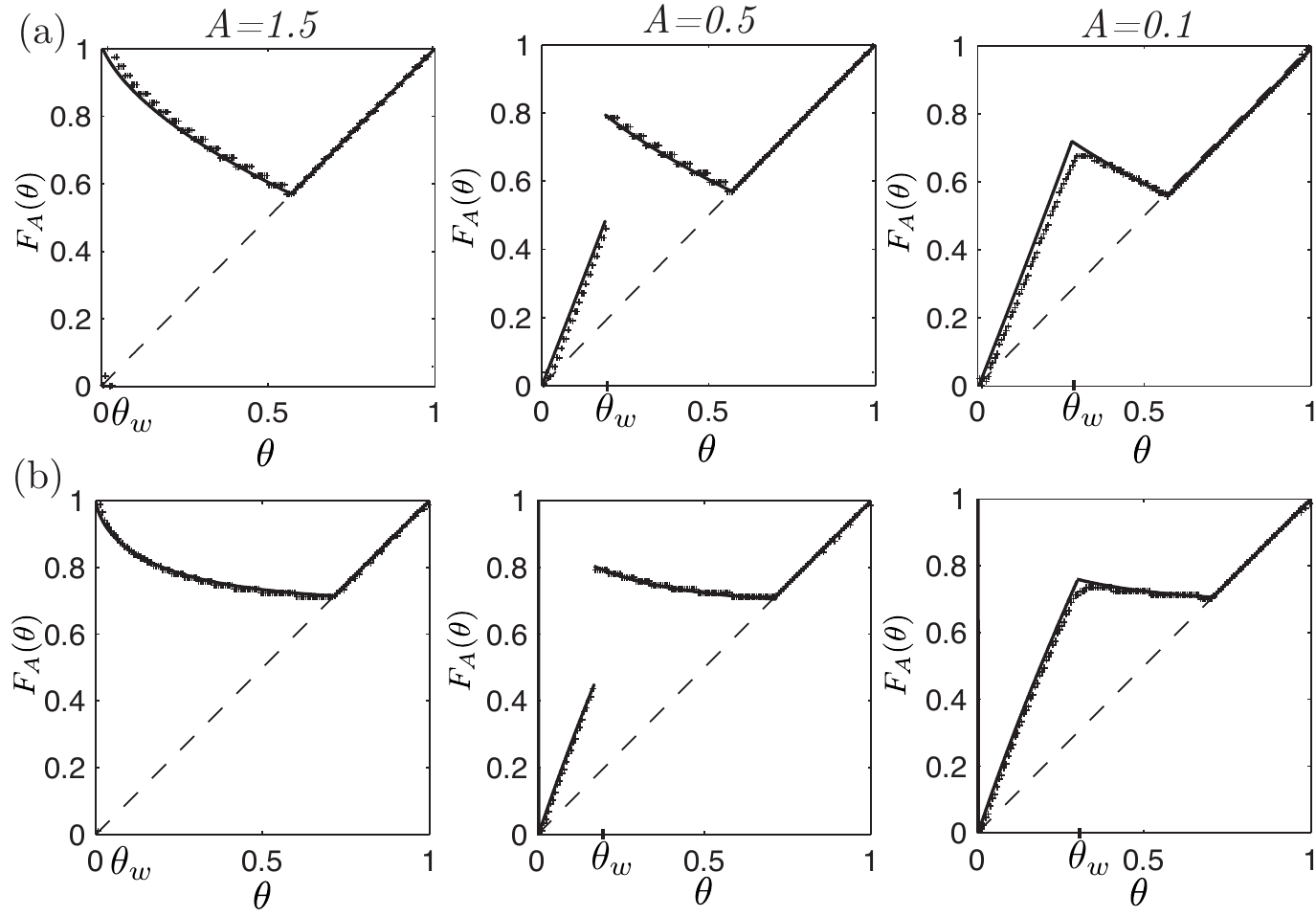}
\caption{Analytically (solid lines) and numerically computed (markers) kick maps $F_A(\o)$ for three values of kick amplitude $A$. (a) Model parameters:  \{$\e=0.01$, $w=1$, $a=0.8$, $b=0$\}. (b) Model parameters:  \{$\e=0.01$, $w=1$, $a=0.4$, $b=0.5$\}.}
\label{fig:computed_map}
\end{center}
\end{figure}
Using the expressions above we can build our kick map by using \eqref{weak_map} and rescaling time by the period of a full cycle. To find the period for a given set of parameters, we use $h_S^{-1}$ and $h_P^{-1}$ to derive the silent and active passage times $T_S$ and $T_P$. We reiterate that the only dependence on a kick's amplitude $A$ is implicitly contained in the expression for $\o_w=h_S^{-1}(y_w(A))$.  As mentioned above, whenever $A>1$, the kicks are strong ($\o_w=0$) and the resulting maps have a fixed shape. On the other hand, as $A$ decreases below 1, the weak kick effect progressively uncovers the left branch of the map ($\o_w \neq0$). 

To verify our derivation of the kick map, we choose 3 prototypical kick amplitudes:  one strong ($A=1.5$) and two weak ($A=0.5$ and $A=0.1$). We then plot the associated maps for two distinct set of model parameters: \{$\e=0.01$, $w=1$, $a=0.8$, $b=0$\} and \{$\e=0.01$, $w=1$, $a=0.4$, $b=0.5$\} in order to better visualize the effects of linear ($b=0$) and saturating exponential ($b>0$) ramps for the $y$ dynamics.

In Fig. \ref{fig:computed_map}, we plot numerically and analytically computed maps, rescaled to the unit circle.  Observe that we get very good agreement between the two and that the main features of weak kicks are captured by our phase reduction model. However, numerically computed maps have fine, plateau-like segments. As argued in \cite{Best:2007p128}, there are as many of these plateaus  as there are spikes in a burst. They appear since a kick can induce bursts with spike counts ranging from one to the number seen in an unperturbed cycle. Since these numbers are integers and we numerically identify phase zero with the last spike of a burst, these plateaus are formed by phases that induce the same number of spikes following a kick. This is not captured by our analytical derivation of the kick map, which is computed from averaged conditions on the fast variable. However, we will see in Sect. \ref{noise} that the presence of noise in the system tends to diminish these plateaus, and in Sect.~5 that key aspects of synchrony and desynchrony for the original ODEs are predicted by our derived kick map.

We close this section with some remarks concerning the generality of features found in the kick map. We use the normal form model (\ref{canon_bautin}) for its  analytical tractability but the general mechanism responsible for bursting, and hence the associated kick map, have characteristics that span across models. For example, for both linear ($b=0$) and saturating exponential ($b>0$) slow dynamics -- which imply distinct interactions between the slow and fast subsystems -- the maps are qualitatively identical.  Specifically, at the end of Sect. \ref{map}, we noted that region of expansion in the left part of the weak kick map is a general feature of elliptic bursters with slow passage effects. The fact that the right branch follows the identity is another general attribute. 

The shape of the middle branch (read left branch for strong kick map) is, however, more model dependent. Notice that both of our strong kick maps (for $A=1.5$) have a left branch that steepens as we move leftwards. This is due to the dependence of $h_P^{-1}$ on $r_P^2(y)$, which grows as $\sqrt{y}$. As a result, for $A$ sufficiently large, $\frac{dF_A(\o)}{d\o}<-1$ for $\o \in (0,\o_c)$ where $\o_c$ is the root of
\beq
\frac{d}{d\o}h_P^{-1}\circ h_S(\o_c)+1=0.
\label{tht_critical}
\eeq
For our parameter set with $b=0$, $\o_c \simeq 0.0968$ while for the second set with $b > 0$, $\o_c \simeq 0.0619$. This curvature shrinks as we decrease the parameter $a$ in (\ref{canon_bautin}) and the silent phase lengthens. 

In general, this branch depends on the slow dynamics and the ratio of silent to spiking times $T_S$/$T_P$, which impose the following constraint:  by construction, a strong kick implies $F_A(0)=1$ and $F_A(T_S)=T_S$. If we were to approximate the $y$ dynamics by constant velocities (as done in \cite{Best:2007p128}), this map branch would be linear and hence contractive, whenever $T_S>T_P$.  In general, as long as the latter is true and the functions $h_S$, $h_P$ (desribing $y$ dynamics) have small total variation, we can expect this branch to be mostly contractive, as for both cases explored above.  We note that we obtain such contraction for the biophysical systems we study in Sect. \ref{GPe}, meant to model GPe neurons in the Parkinsonian state; \cite{Best:2007p128} draws an interesting contrast with cases having $T_S<T_P$.
The generality of these features motivate the analysis of dynamics induced by the kick map as we show in the next section.

In light of these remarks and in the interest of clarity, we use the maps computed above for the parameter set with $b=0$ to carry out our analysis for the rest of the paper. We stress that the arguments that will follow hold true for other parameter sets of system \eqref{canon_bautin} and any elliptical bursting system having the features described above.

\section{Dynamics of the kick map }
\label{dynamics}

Now that we have an understanding of an elliptic burster's response to input kicks of various strengths $A$, we turn to the other input parameter of relevance --  the period between these kicks $\tau$ -- and study the iterated dynamics of the map for various combinations of $A$ and $\tau$.  

\subsection{Iterative framework}

\begin{figure}[h]
\begin{center}
\includegraphics[width=250pt]{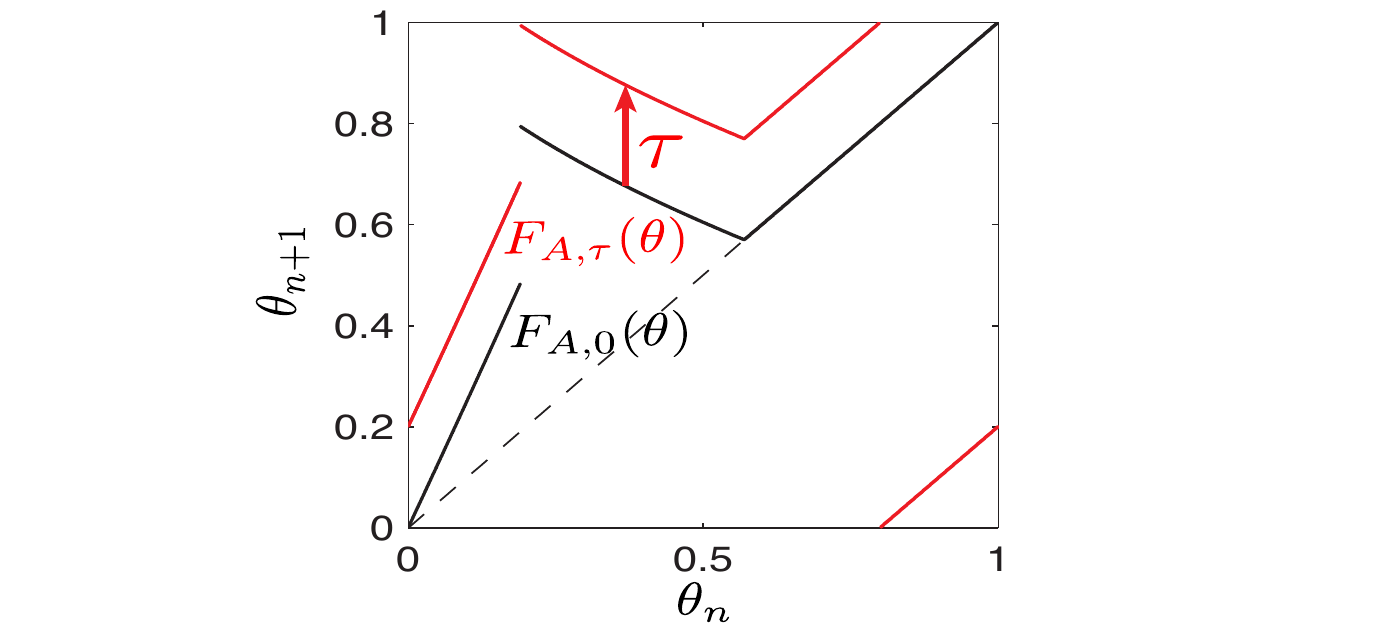}
\caption{Effect of kick period $\tau (\mod 1)$ on the kick map. Original kick map $F_A(\o)$ in black and $\tau$ induced kick map $F_{A,\tau}(\o)$ in red. $\tau$ simply translates the map vertically.}
\label{fig:translate}
\end{center}
\end{figure}

We now use the kick map to build an iterative dynamical system capturing the evolution of cells subject to periodic stimulation.  Let 
\begin{equation}
F_{A,\tau}(\o) \equiv F_A(\o)+\tau \; (\textrm{mod} \;1)
\label{kick_map}
\end{equation}
which returns the phase of a kicked cell right before the next kick, $\tau$ time units later. The parameter $\tau$ translates the map vertically as illustrated in Fig. \ref{fig:translate}. For a chosen pair $(A,\tau)$ and some initial phase $\o_0$, an orbit is defined by $\o_{n+1}=F_{A,\tau}(\o_{n-1})=F^{n}_{A,\tau}(\o_0)$.

We note that, in contrast with the map of~\cite{Best:2007p128}, our kick map does not become rescaled as $\tau$ is varied, but is rather translated around the circle (as in related studies~\cite{Belair:1986p8707,Glass:1983p9985}).
Another difference is the presence of sustained expansion even though $T_S>T_P$ (see above), due to accelerated slow passage effects induced by weak kicks; we will show that this expansion leads to positive Lyapunov exponents for certain values of $\tau$. This phenomenon has been exploited in simpler dynamical systems with slow passage through a supercritical Hopf bifurcation in the context of chaos control \cite{Perc:2006p5100}. 

Equation~\eqref{kick_map} assumes that a kick acts on cells located on unperturbed trajectories.  As also noted above, the separation of timescales for elliptic bursters implies very fast attraction back to steady state trajectories following a kick, so that this assumption is generally valid.  Moreover, the iterated dynamics for small values of $\tau$ are relevant in any case, as they can be accessed by longer, equivalent kick periods modulo one.

We next use Eqn.~\eqref{kick_map} to study the response of a population of identical elliptic bursters with different initial conditions to a common, pulsatile signal.  
For example, globally stable fixed points or periodic orbits represent phase locking regimes towards which the long term behavior of any cell will converge. On the other hand, maps that yield sensitivity to initial conditions and complex orbits are representative of desynchronizing inputs, when delivered to a population. 
To further explore population behavior, we need to define a metric by which we quantify synchrony of phase points on $S^1$. We call this our ``synchrony measure" and describe it next.

\subsection{Assessing synchrony}

There are many ways one can quantify how closely $N$ points are distributed on $S^1$. Two natural choices are the binned entropy $H$ and order parameter $R$ (also known as vector strength) :
\begin{equation*}
\begin{split}
H(\o_1,...,\o_N) &= \frac{1}{\log(1/N)}\sum_{j=1}^Np_j \log(p_j)\\
R(\o_1,...,\o_N) &=|\frac{1}{N} \sum_{j=1}^N e^{i2\pi \o_j}|.
\end{split}
\end{equation*}
For $H$, we divide $S^1$ into $N$ equal length subintervals, or bins, and take $p_j$ to be the number of phases in bin $j$ over $N$ (using the convention $0 \log 0 = 0$). $H$ takes its maximal value one when there is a phase spread into each bin and its minimum value zero when all are concentrated into a single bin. On the other hand, $R$ takes its minimum value zero when phases are evenly distributed, and one when they are all equal. 

Each measure has strengths and weaknesses as a metric of synchrony.  For example, $R$ can be zero if the phases are split into two equal, antipodal groups on $S^1$ -- this is hardly an asynchronous state. For a large $N$, $H$ can take relatively large values even if cells are distributed in close by bins.  By taking  
\beqn
\begin{split}
W(\o_1,...,\o_N) &= \frac12[R(\o_1,...,\o_N)+(1-H(\o_1,...,\o_N))].
\end{split}
\eeqn
to be the average of $R$ and $1-H$, we can be assured that low measures of $W$ correspond to cases where cells are well distributed across bins and that these bins are broadly spread around $S^1$.  Specifically, we assess synchronization properties of a given map by taking $N$ cells $\{\o_n\}_{1 \leq n \leq N}$ with some initial distribution on $S^1$, pushing these states forward through $m$ iterates, and computing $\bar W$ as the average over the last $k$ out of these $m$ iterates:  
\beq
\bar{W}=\frac{1}{k}\sum_{i=m-k}^m W(F^i(\o_1),...,F^i(\o_N)).
\label{synch_avr}
\eeq
Throughout the paper, we take $m \geq 100$ and $k=20$, having found empirically that values change little with larger values of either.

\subsection{High period orbits and positive Lyapunov exponents}

We next investigate how iterations of our kick maps act on a population of cells for the three prototypical cases of strong ($A=1.5$) and weak kick maps ($A=$0.5, 0.1), plotted in panel (a) of Fig. \ref{fig:computed_map}.
We note that, while smooth maps on the circle are well characterized~\cite{Strogatz:1994p9396, Katok:1995p9450}, the discontinuities in our map introduce a number of distinct phenomena -- such as border collision bifurcations with period adding at all orders, and ``sharp" transitions to chaos. There is an ongoing effort to build a theory to better understand such systems \cite{DiBernardo:2008p9459,Boyarsky:1997p6982,Hogan:2007p3847,Jain:2003p3862}.

\begin{figure}[t!]
\begin{center}
\includegraphics{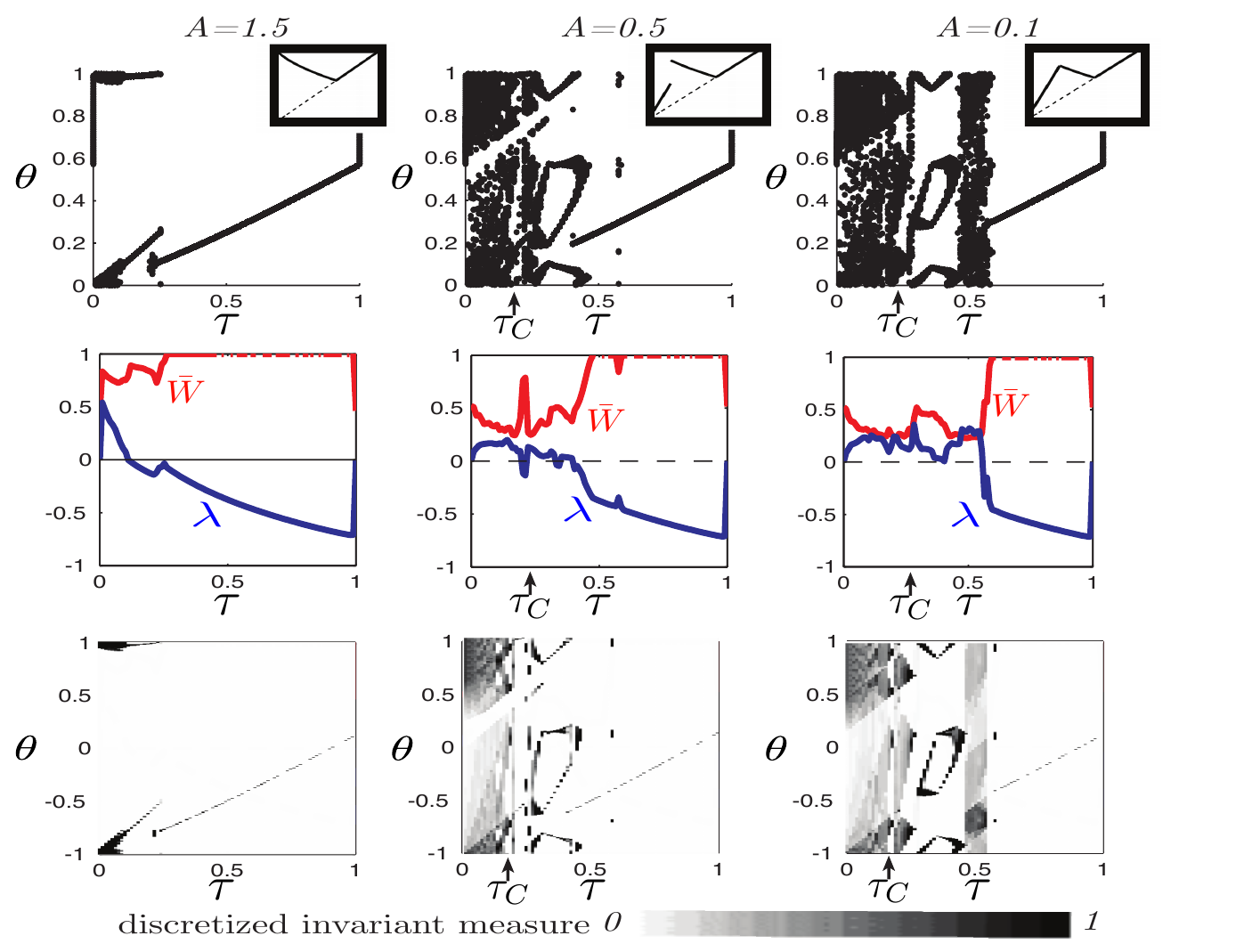}
\caption{Orbit diagrams and related measures for kick maps (Eqn.~(\ref{weak_map})). Top to bottom:  orbit diagrams of 100 cells after 150 iterations;  synchrony measure $\bar{W}$ (red) and averaged Lyapunov exponent $\l$ (blue);  invariant measure approximates. Left to right:  strong kick with $A=1.5$;  weak kick with $A=0.5$;  weak kick with $A=0.1$. Marked values for $\tau_C$ below which Lemma \ref{lemma_chaos} applies.}
\label{fig:bif_diag}
\end{center}
\end{figure}

We first plot orbit diagrams with respect to the parameter $\tau$ for each of the three maps at hand (Fig. \ref{fig:bif_diag}). Specifically, we select 100 cells uniformly distributed on $S^1$; for various $\tau \in [0,1]$, we compute the positions of these cells after 150 iterates of $F_{A,\tau}(\o)$ and ``vertically" plot the result.
Directly below these orbit diagrams, we plot the synchrony measure $\bar{W}$ of these end states, together with numerically computed Lyapunov exponents $\lambda$ for each $\tau$, averaged over all trajectories.
Finally we compute approximations of invariant measures for each $\tau$, using a variation of Ulam's method developed in \cite{Ding:2001p4923} (this produces a discretized approximation of fixed densities for a map's Perron-Frobenius operator). The results are plotted in the bottom panels of Fig. \ref{fig:bif_diag}.

We begin by describing results for the strong kick map ($A=1.5$).  Recall that the leftmost part of the strong kick map ,when $\o \in (0,\o_c)$ ($\o_c \simeq$0.0968, derived in Eq.~\eqref{tht_critical})  has a derivative greater than one in absolute value. The derivative is less then or equal to one in absolute value everywhere else. It is easy to see that the map has a fixed point for any $\tau$. If the map intersects the identity at $\o \in (0,\o_c)$, the fixed point is unstable and we find that complex dynamics emerge. We can witness this by looking at the orbit diagram of this map for small $\tau$. Nevertheless, we see that cells tend to cluster in small regions and hence are relatively synchronized. For any other $\tau \neq 0$, we have stable fixed points, implying phase locking of cells to the input kicks. Note that for all maps, $\tau=0$ implies a continuum of neutrally stable fixed points, as the right branches of the maps align perfectly with the identity.

For the weak kick maps ($A=0.5, \; 0.1$), some values of $\tau$ yield stable, discrete attractors as well:  fixed points and periodic orbits.  More interesting are values of $\tau$ which produce thick, chaotic attractor-like objects. These are associated with what appears to be locally absolutely continuous invariant measures and positive Lyapunov exponents. It is in these regimes that our synchrony measures show a greater spread of cells. This is indicative of chaos, the presence of which is consistent with expansive regions of weak kick maps. Interestingly, these regions appear wether a map has a gap or not, as we show below. 

Rigorously assessing the presence of chaos is, however, not a simple task.  This was accomplished for related piecewise smooth maps~\cite{Bressloff:1990p8661,Coombes:1996p8613,Keener:1980p4733,BUB:1995p8616} where various definitions of chaos were used, depending on context. For example, Keener \cite{Keener:1980p4733} showed that piecewise, surjective and non-decreasing maps with overlap have rotation numbers spanning a non-empty interval, an indication of chaos for circle maps. Unfortunately, our weak kick maps sometimes fail to have surjectivity (e.g. A=0.5) and always fail to be non-decreasing. One can also try to define trapping regions and a family of intervals for which interval images cover at least one other interval and the image of at least one interval covers at least two others.  This constructs an shift on a space of sequences, which can characterize a chaotic system. However, building such an interval family for the weak kick map proves to be quite complex as intervals get flipped by decreasing parts of the map and severed by discontinuities.  Thus, we do not aim at a complete characterization of the complex dynamics produced by our weak kick maps.  However, we next show that, for some $\tau$ values, there is a positive lower bound on (sustained) Lyapunov exponents of any trajectory.

In Fig. \ref{fig:bif_diag}, we see that some values of $\tau$ ($0 < \tau \lesssim0.2$) induce dynamics that appear chaotic for both of our weak kick maps.  Thus motivated, we state the following lemma:

\begin{lemma}
Consider a piecewise defined map $F(\o)$ on the circle which is smooth on three non-intersecting intervals $I_1$, $I_2$ and $I_3$ with $\bar{I}_1 \bigcup \bar{I}_2 \bigcup \bar{I}_3 = S^1$. Suppose $|\frac{dF}{d\o}|_{I_1}| \geq a >1$, $|\frac{dF}{d\o}|_{I_2}| \geq b >0 $ -- such that $\ln(a) > |\ln(b)|$ -- and additionally that $\frac{dF}{d\o}|_{I_3}=1$. Then if  $F(I_2) \subset I_3$, $F(I_3) \bigcap I_2 = \emptyset$ and $F(I_3) \neq I_3$, the Lyapunov exponent associated with the orbit of almost any initial condition $\o_0 \in S^1$, if well defined, will be strictly greater than zero ($\l(\o_0)>0$).
\label{lemma_chaos}
 \end{lemma}
 \textbf{Proof:}
 Given $\o_0$ and its forward orbit $\{\o_n\}_{n=0,1...}$, the local Lyapunov exponent can be written as $ \l(\o_0)=\lim_{N \to \infty} \frac1N \sum_{n=0}^N \ln|\frac d{d\o} F(\o_n)|$. The derivative is defined everywhere in $S^1$ except the 4 border points of the intervals $I_1,I_2,I_3$, which is obviously a measure zero set. The condition $F(I_2) \subset I_3$  imply that any point in $I_2$ is sent to $I_3$. Since $F$ is smooth on $I_3$ and $\frac{dF}{d\o}|_{I_3}=1$, the condition $F(I_3) \neq I_3$ implies that any point in $I_3$ must eventually exit it. Let $C=\max_\o \{ n= \min_m \{m \, | \,F^m(\o) \notin I_3 \}\, | \, \o \in I_3\}$. Then any element of $I_3$ stays in $I_3$ at most $C$ iterates. Furthermore, $F(I_3) \bigcap I_2 = \emptyset$ implies that elements of $I_3$ are eventually sent to $I_1$.
When an orbit point visits $I_1$, it contributes at least $\ln(a)>0$ to the sum in $\l(\o_0)$, at least $\ln(b)$ (possibly$<0$) if it visits $I_2$ and 0 when it passes by $I_3$. By tracking which intervals an orbit visits, any admissible subsequence featuring $I_2$ must contain $I_2 \to I_3$. The sequence that contributes the least to $\l(\o_0)$ is therefore $I_2 \to I_3 \to ... \to I_3 \to I_1$ where $I_3$ is repeated $C$ times. It follows that for almost every $\o_0$, $\l(\o_0) \geq \frac1{2+C}(\ln(a)+\ln(b)+0+...+0) >0$. $\square$
 \bigskip

\begin{figure}[h]
\begin{center}
\includegraphics{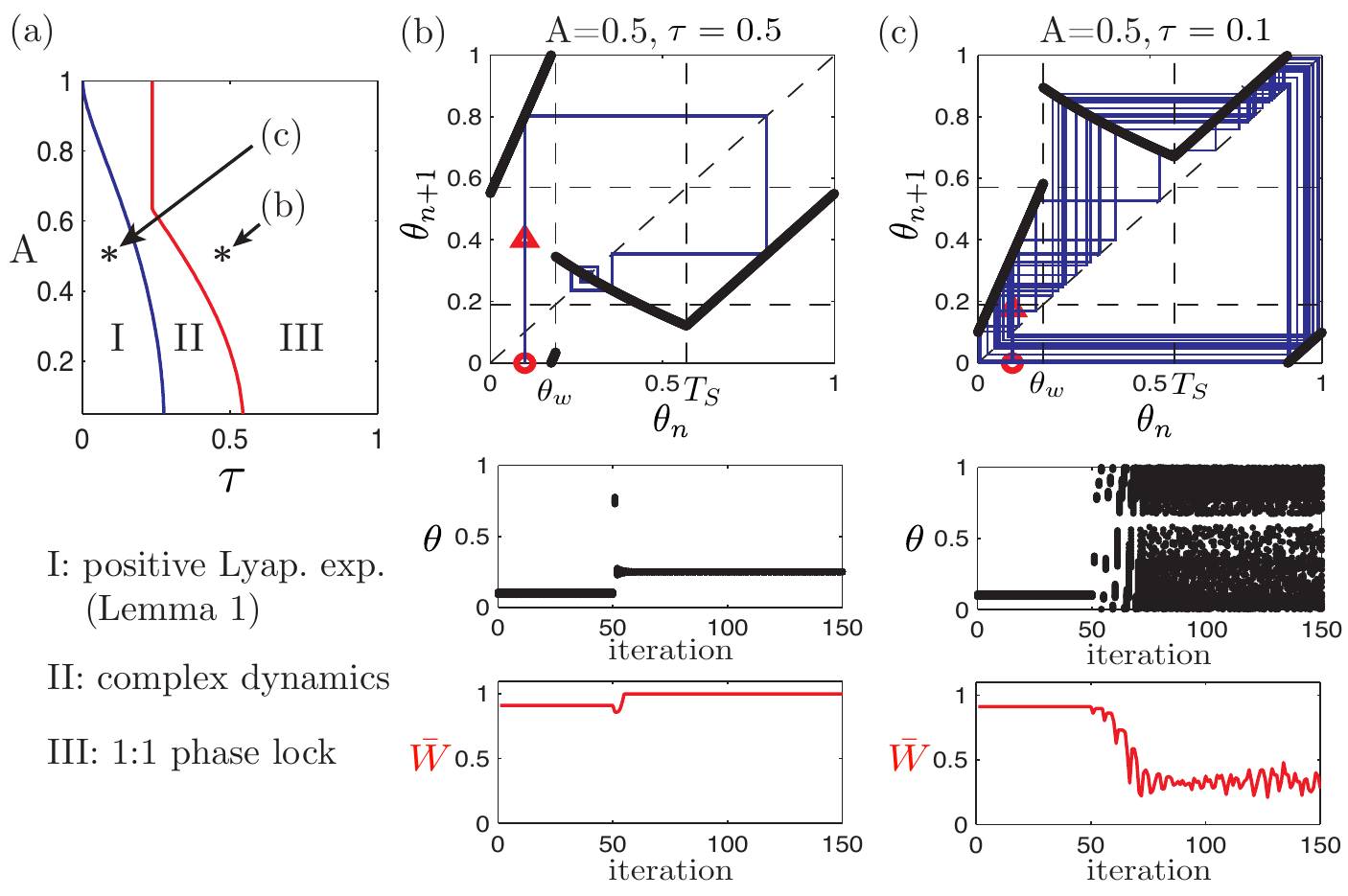}
\caption{(a) Three regions of map dynamics in $(A,\tau)$-space (see text for details). (b) \& (c):  two sample evolution trajectories for the kick map with $A=0.5$. Top to bottom:  cobweb diagrams and marked intervals $(0,\o_w)$, $(\o_w,T_S)$ and $(T_S,1)$ for which Lemma \ref{lemma_chaos} applies;  sample trajectories with 50 initial condition randomly chosen close to $\o=0.1$ (iterations 50 and above are from kick map);  synchrony measure $\bar{W}$. (b) $\tau=0.5$, cells synchronize;  (c) $\tau=0.1$, cells desynchronize.  }
\label{fig:trajectories}
\end{center}
\end{figure}

We note that weaker conditions could be stated under which a positive Lyapunov exponent results, but the above are sufficient for the map at hand, as we now show. 
In particular, we show that Lemma \ref{lemma_chaos} can be applied to the weak kick maps for certain values of $\tau$. Clearly, the intervals $(0,\o_w)$, $(\o_w, T_S)$ and $(T_S,1)$ will play the roles of $I_1$,$I_2$ and $I_3$. For our prototypical weak kick maps, we have $|\frac {dF_{A,\tau}}{d\o}|_{(0,\o_w)}| \geq 2.7$, $|\frac {dF_{A,\tau}}{d\o}|_{(\o_w, T_S)}| \geq 0.65$ and $\frac {dF_{A,\tau}}{d\o}|_{(T_S,1)} =1$ which fulfills the derivative criteria.  It remains to show that the intersection requirements for interval images are met; we take a graphical approach which leads to conditions on $\tau$. In Fig. \ref{fig:trajectories}, panels (b) and (c) show the map for $A=0.5$ along with the marked intervals for two distinct values of $\tau$ and cobweb diagrams of sample trajectories. One can verify that the Lemma's assumptions are respected as long as the leftmost tip of the middle branch stays smaller than one and the rightmost tip is greater than the identity (panel (c)). This happens when $0<\tau<\tau_C$ where 
\begin{equation}
\tau_C=1-\lim_{\o \to \o_w^+}F_{A,0}(\o).
\label{tau_c}
\end{equation}
 Although we do not explicitly graph it, the same argument holds for $A \in (0,1)$. We analytically compute these upper bounds for $\tau$ and get $\tau_C\simeq0.205$ for $A=.05$ and $\tau_C \simeq0.26$ for $A=0.1$. As expected, these are close to $0.2$, the rough higher bound we predicted earlier from Fig.~\ref{fig:bif_diag}. In addition to positive Lyapunov exponents, $\tau<\tau_C$ imposes a cyclic structure where trajectories visit a large portion of all three intervals in finite time, which is necessary for the population to become widely distributed around $S^1$.

For $\tau > \tau_C$, various dynamics can be observed. Figure \ref{fig:bif_diag} shows that positive Lyapunov exponents can still be found sporadically but the regions on which the trajectories accumulate can be considerably smaller. Periodic orbits of various periods are also present; in particular, high $\tau$ values seem to be associated with stable fixed points. This motivates a separation of the $(A,\tau)$--space into three regions:  I, where Lemma \ref{lemma_chaos} applies;  II, where dynamics are complex and transitions between what appears to be periodic orbits and smaller chaotic regions can be seen;  III, stable fixed points (1:1 phase locking). Panel (a) in Fig. \ref{fig:trajectories} shows these regions. We stress that further analysis of region II might yield additional structure but the complexity of our map (circular domain, increasing and decreasing parts, transitions from gaps to discontinuity points, etc.) renders a complete analysis outside of this article's scope.

From what was described above, we see that region I is given by 
$$\textrm{I}: \, \{(A,\tau) | 0<\tau<\tau_C(A)=1-\lim_{\o \to \o_w(A)^+}F_{A,0}(\o)\}$$
where we have written the expression for $\tau_C$ to highlight its dependence on the kick amplitude $A$. As $A$ increases, $\o_w$ decreases and consequently, $\tau_C$ as well. Region I vanishes altogether when the kick transitions from weak to strong -- at $A=1$, where $\tau_C=0$.  

To define the boundary between regions II and III, we need to derive a condition under which our maps have a stable fixed point. We begin by inquiring about when the middle branch of the map intersects the identity.  This happens when $\tau_C+\o_w<\tau<1$. The stability of the resulting fixed point depends on the derivative of the branch at the intersection point. 

Again, as $A$ increases and both $\o_w$ and $\tau_C$ decrease, more of the middle branch of the map is exposed. When $\tau=\tau_C+\o_w$, the leftmost tip of this branch intersects the identity. If $A$ is sufficiently large, the resulting fixed point will be unstable; recall that for the strong kick map ($A>1$), $|\frac{dF}{d\o}|>1$ when $\o \in (0,\o_c)$. It follows that if $\o_w<\o_c$, the first fixed points to appear as $\tau$ increases are unstable. We include unstable fixed points in region II and get the following definitions
\beqn
\begin{split}
\textrm{II}:& \{(A,\tau) | \tau_C(A)<\tau< \tau_C(A)+\max\{\o_w(A),\o_c\}\}\\
\textrm{III}:& \{(A,\tau) | \tau_C(A)+\max\{\o_w(A),\o_c\}<\tau< 1\}.
\end{split}
\eeqn

What can be taken from our analysis thus far is that to desynchronize cells with a $\tau$-periodic input, the best strategy appears to be weak kicks with $0<\tau<\tau_C$ (region I). To illustrate this, panels (b) and (c) of Fig. \ref{fig:trajectories} show sample trajectories for 50 cells under the action of the weak kick map with $A=0.5$ in both synchronizing (region III) and desynchronizing (region I) regimes.  To mimic a synchronous state, our initial phases are drawn from a uniform distribution on an interval of width .01, around $\o=1/10$. The first 50 iterates are taken with respect to the identity map, reflecting a preliminary period with no pulsatile inputs. The next iterates are taken with respect to the weak kick map with $\tau=0.5$ (region III, panel (b)) or $\tau=0.1$ (region I, panel (c)). 

As expected, for $\tau=0.5$, every cell attracts to a single fixed point and the synchrony measure is maximized. However, for $\tau=0.1$, the cells quickly become widely distributed. This was foreseeable from the shape of the computed invariant measures at this value of $\tau$. Importantly, the synchrony measure drops substantially. 

We close this section with an important remark concerning circle maps and phase locking. Here, our kick map can be seen as a $A$-perturbation of a $\tau$-rotation. Typically, smooth perturbations of rotations admit Arnold Tongues:  well separated wedge like regions where phase locking with rational rotation number occurs (see e.g. \cite{Wiggins:2003p12781}, Sect. 21.6). In our case, the loss of smoothness introduces sustained region where chaotic dynamics prevail, hereby showing surprising effects of discontinuous circle maps.  \\

\section{Effects of noise on kick map and synchrony}
\label{noise}

Up to this point, we developed a phase reduction framework to analyze the dynamics of periodically forced bursters.  We now ask:  does the phase reduction remain valid in the presence of noise?  If so, what qualitative changes in the map occur, and what are the consequences for entrainment of bursters?  In this section, we develop answers through a blend of numerical results and analytical approximations.

The effects of stochastic perturbation on elliptic bursters have been previously studied in various contexts~\cite{Baer:1989p4834,Su:2004p2369, Kuske:2002p11725}. The unifying theme is the effect of noise on slow passage through a Hopf bifurcation.
Here, we build on these results to better understand the role of noise on burst responses to periodic pulsatile inputs.  We will comment further on prior results as we progress.
  
To this end, we introduce a stochastic perturbation in the fast variables of the normal form model, so that Eqn.~\eqref{canon_bautin} becomes  
\begin{equation}
\begin{split}
\dot{z}&=(y+iw)z+2z|z|^2-z|z|^4+I(t)+ \n \xi(t) \\
\dot{y}&=\e(a-|z|^2-by).
\end{split}
\label{canon_noisy}
\end{equation}
where $\n\geq 0$ is a small, real, noise strength parameter and $\xi(t)$ is a time-periodic train of discrete small kicks with normally distributed amplitudes. That is, $\xi(t)=\sum_i \xi_i \d(t-i\D t)$, where the $\xi_i$ are i.i.d. as $N(0,\sqrt{\D t})$ and the timestep $\D t$ controls the temporal resolution of the noisy perturbation. We take $\D t=0.05$, the numerical solver's maximal timestep, so that $\xi(t)$ approximates white noise (see Sect. \ref{numerics} for more on numerical methods). Notice that the noise term only acts on the real part of the fast variable $z$, which mimics a cell's voltage variable.   

In~\cite{Baer:1989p4834}, the authors show that such noise terms diminish slow passage effects through supercritical Hopf points (i.e., causing cells to jump to the spiking state closer to $y_H$); in~\cite{Su:2004p2369}, a similar effects was found for elliptic bursters \cite{Su:2004p2369, Kuske:2002p11725}. As we discuss below, there are cases where a phase reduction can still be defined in the presence of such noise, with an interesting and tractable impact on the kick map's shape.

\subsection{Effects of noise on the phase reduction of elliptic bursters}

Our previously derived phase reduction relied on an important assumption:  periodicity of the burst cycle. As discussed in Sect. \ref{geometry}, elliptic bursters do not necessarily have periodic solutions, but rather a metastability property which guarantees the constancy of cycle's duration $T$, up to $O(\e)$.  This regularity is what enables our phase reduction.   


When noise is added to the fast subsystem, either regular or highly variable burst durations can result, depending on noise strength and type~\cite{Baer:1989p4834,Su:2004p2369, Kuske:2002p11725}.  Below, we will show that there is a wide range of noise strengths that significantly impact the underlying dynamics, but maintain regular burst durations.  Specifically, for system \eqref{canon_noisy} ($\e=0.01$, $w=1$, $a=0.8$, $b=0$) we compute the coefficient of variation (CV) of these durations as described in Sect.~\ref{redux}, for noise strengths $\n$ ranging from $10^{-17}$ to $10^{-1}$. 

Panel (c) of Fig. \ref{fig:noisy_stats} shows our findings. Although the mean period $\la T \ra$ decreases with increasing noise strength, the $CV$ remains low -- below $10^{-2}$ -- for the range of noise strengths $\n \leq 10^{-3}$ (recall that $CV \approx 10^{-3}$ for the noiseless case).
 In other words, for a wide range of noise strengths, random forcing does not introduce substantial variability to the burst period. This is consistent with results from~\cite{Su:2004p2369} where the authors show in a closely related setting that the silent and spiking times $T_S$ and $T_P$ are related to the log of the noise amplitude.  

Despite the fact that they preserve regular burst periods, noise strengths $\n \leq 10^{-3}$ have a strong impact on the slow passage effect.  
The integral condition~(\ref{integral_condition}) reflects cancellation of attraction to $S$ by and repulsion away from it; noise limits the extent of attraction and therefore the duration required for repulsion.  This is illustrated in panel (a) of Fig.~\ref{fig:noisy_stats}, which shows that the jump up point from silent to spiking regimes decreases as $\n$ increases within $[0,10^{-3}]$.  In this regime, the averaged $y$ dynamics are relatively unaffected by this stochastic forcing, so that $CV$ remains low.   Once $\n = 10^{-2}$, there is no slow passage and the solutions jump up to spiking as soon as $y$ crosses $y_H$. For $\n>10^{-2}$, the stochastic kicks have accumulated effects comparable to our forcing kicks and we see solutions randomly jumping into spiking before $y$ reaches $y_H$, which explains the increasing $CV$ for this range of noise strengths~\cite{Su:2004p2369}. 

We next pursue phase reduction to an approximate, deterministic circle map for the low $CV$ cases ($\n \in [0,10^{-3}]$).

\subsection{Effect of noise on the kick map}
\label{noisy_DD}

For strong kicks -- where responses are not determined by slow passage effects  -- adding noise does not considerably change the shape of the kick map; the chief effect is that the region where $| \frac{dF_A}{d\o}|>1$ shrinks (not shown). This similarity was expected since a strong kick instantaneously translates a bursting solution to the active phase, where noise has little effect since fast dynamics follow large amplitude trajectories.

\begin{figure}[h]
\begin{center}
\includegraphics{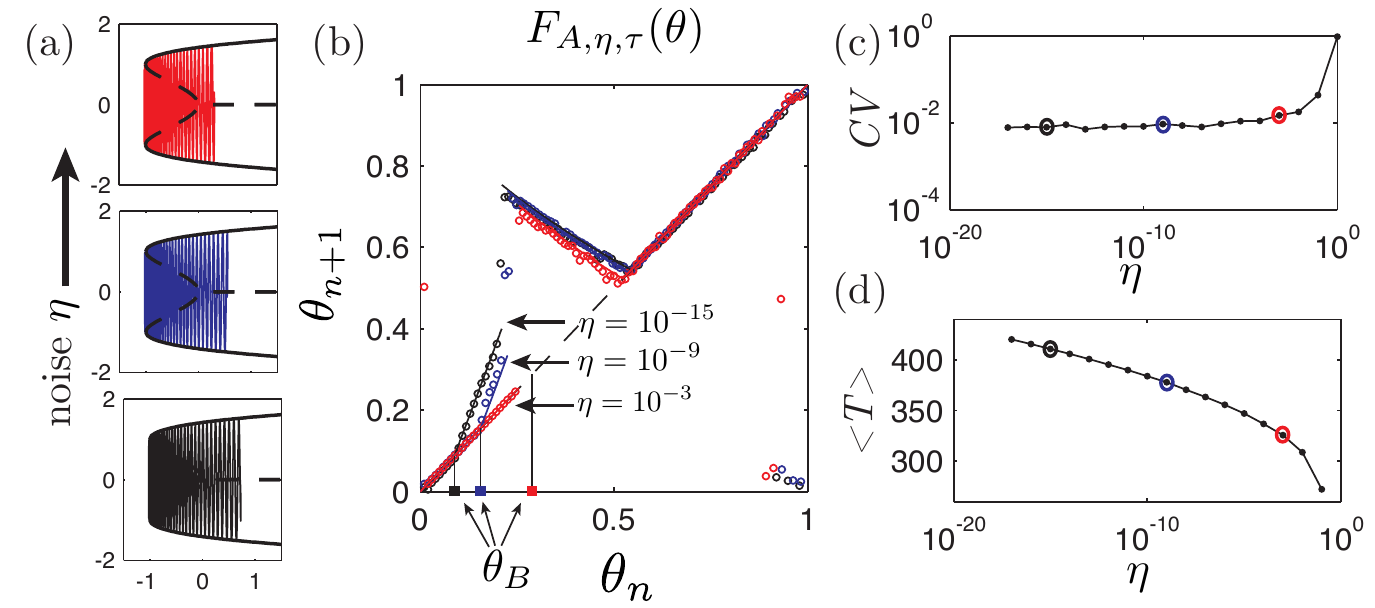}
\caption{(a) Slow passage effect shortens with increasing noise. (b) Numerically and analytically computed kick maps (
from Eqn. \eqref{noisy_map}) for noise strengths $\n=10^{-15}$, $\n=10^{-9}$ and $\n=10^{-3}$ (all with kick amplitude $A=0.5$). (c) $CV$ of burst cycle period (computed using 150 cycles) against $\n$ using a log-log scale. (d) Average burst cycle period $\la T \ra$ (same sample as in (c)) against $\n$ using a log-linear scale.}
\label{fig:noisy_stats}
\end{center}
\end{figure}

For the remainder of this section, we concentrate on the considerable changes noise induces for weak kick maps. In panel (b) of Fig.~\ref{fig:noisy_stats} we plot numerically computed maps for $A=0.5$ and $\n =$ $10^{-15}$, $10^{-9}$, $10^{-3}$. Each marker represents the average response of a given initial phase to 10 distinct realizations of stochastic forcing. In what follows, we derive approximations for these maps, plotted in solid lines on the same figure. In order to proceed, we discuss key differences among numerical maps of stochastic bursters with varying noise strength.

We can see that the maps are qualitatively similar in almost all aspects except the left branch, which collapses on the identity as noise increases. Recall that the expansion in this branch is due to altered slow passages due to weak perturbations on the stable part of the silent branch $S$. As discussed above, noise shortens slow passages and one might expect the steepness of the expansive part to decrease as noise increases. Surprisingly, it is not this steepness that changes but rather the onset of the expansive ramp ($\o_B$ defined below) which varies. To explain this phenomenon, we must turn to the concept of buffer points for delayed bifurcations.

When deriving the analytical kick map \eqref{weak_map}, we relied on the integral condition \eqref{integral_condition} which dictates that the further a solution starts (at $y_i$) from the Hopf point $y_H$, the longer the slow passage will be. This is true for $y_i$ in a range before $y_H$, up to a point beyond which this relation fails. In fact, given some $y$ dynamics, one observes that the length of the slow passage ($y_j-y_H$) saturates to a constant value for any initial $y_i$ far enough from $y_H$. This is called the \textit{maximal delay}. The closest point to the ``static" bifurcation point ($y_H$) to generate a maximal delay is called the {\it buffer point} of a delayed bifurcation~\cite{Diener:1991p11863}. The techniques used in~\cite{Baer:1989p4834,Su:2004p2369, Kuske:2002p11725,Baer:2008p11521} to analyze deterministic and stochastic slow passages assume that a solution always remains closer to the bifurcation than the associated buffer point.  This is also the case for the bursters we study in the absence of noise, where the buffer point is to the left of $y_{SN}$.  Thus, it does not affect the dynamics, enabling us to use condition~\eqref{integral_condition}. 

For noisy bursters, the numerically derived weak kick maps in Fig.~\ref{fig:noisy_stats} show that the buffer point can lie to the {\it right} of $y_{SN}$, therefore exerting an important effect. If a weak kick is delivered when $y_i$ is to the left of the buffer point, the trajectory retains no memory of the kick and the phase of the cell is left unchanged by the kick.  If it is delivered when $y_i$ is to the right of the buffer point, the kick will shorten the slow passage, as for the noiseless case. 

To our knowledge, there are no general results in the literature deriving buffer points for noisy delay bifurcations. In what follows, we derive an approximation for them, in the context of weakly kicked trajectories, that holds for Eqn.~\eqref{canon_noisy}. We work in the small $\Delta t$ limit, so that the stochastic forcing is white noise.

\subsubsection{Buffer points and stochastic slow passage through the Hopf point}

For a given noise strength $\n$, define $y_B$ to be the buffer point for the delayed Hopf bifurcation in Eqn. \eqref{canon_noisy}. That is, $y_B$ is the smallest initial value $y_i$ such that a weak kick delivered at $y_i$ will induce a change in the jump up point $y_j$. Naturally, $y_B<y_H$ and judging by the shape of our weak kick maps in Fig. \ref{fig:noisy_stats}, we can expect $y_{SN}<y_B$. We now derive a probabilistic criterion that gives an accurate approximation for $y_B$.

We begin by assuming that $\Delta t$ is sufficiently small in Eqn. \eqref{canon_noisy} so that the noise term  $\n \xi(t)$ can be approximated by the white noise term  $\n dW(t)/dt$, 
where $W(t)$ is a real valued Wiener process. We are interested in solutions following the silent branch $S$ when $y \in (y_{SN},y_H)$. As done in Sect. \ref{canon_map}, we make the singular limit assumption that $z$ travels on $S$ and therefore that $|z|=0$.  This enables us to decouple the $y$ dynamics and get an expression for $y(t)=y_i+g(\e t)$, precisely as in the deterministic case (see Eqns. \eqref{fret0},\eqref{fret}).  These assumptions hold in the limit that $\e$ and $\n \rightarrow 0$, and we will show that they give good approximations for the parameters used here ($\eps =0.01$, $\eta \le 10^{-3}$).

Since we are interested in solutions near $S$, which is composed of fixed points of the fast subsystem, we linearize the fast dynamics about $z=0$. We write the resulting equation in real coordinates $x=(x_1,x_2)^T$ where $z=x_1+ix_2$:
\beq
dx=J(t)xdt+BdW(t)
\label{canon_linear}
\eeq
where
\beqq
J(t)=\left(\begin{array}{cc}y(t) & -w \\w & y(t)\end{array}\right), \quad B=\left(\begin{array}{cc}\n & 0 \\0 & 0\end{array}\right)
\eeqq
and $W(t)$ is now a two dimensional real valued Wiener process. Notice that $J(t)$ commutes with itself ($J(t)J(s)=J(s)J(t)$ $\forall$ $t,s$) which enables us to write the noise-free ($\eta=0$) solution of \eqref{canon_linear} as
\beqq
x(t)=e^{\int_0^tJ(s)ds}x_i
\eeqq
where $x(0)=x_i$. Using this property, Ito's formula yields explicit expressions for the mean $\mu(t)$ and covariance matrix $\Sigma(t)$ of the time dependent probability distribution for $x$ governed by Eqn. \eqref{canon_linear}:
\beqr
\mu(t)&=& e^{\int_0^tJ(s)ds}\mu_i \label{plo1}\\
\Sigma(t)&=& e^{\int_0^tJ(s)ds} \Sigma_i e^{\int_0^tJ(s)^Tds}+\int_0^t dt' e^{\int_{t'}^tJ(s)ds}BB^Te^{\int_{t'}^tJ(s)^Tds}\label{plo2}
\eeqr
where $\mu(0)=\mu_i$ and $\Sigma(0)=\Sigma_i$. See \cite{citeulike:1400625}, Chap. 4 for details of this derivation.
In particular, we will suppose that we have an initial distribution for $x$ that is {Gaussian}; then, the distribution of $x$ at any time $t$ is fully determined by Eqns. \eqref{plo1} and \eqref{plo2}, as a Gaussian distribution with mean $\mu(t)$ and $\Sigma(t)$.  We introduce the following notation $\tilde p(\cdot)$ for this distribution, which makes clear the dependence on the initial condition for $y(0)=y_i$ as well as the initial distribution of $x$ and the elapsed time $t$: 
\beq 
x \sim \tilde p(x(t) | y_i, \mu_i, \Sigma_i) \; \;. 
\label{e.gen}
\eeq

We next study the distribution of our (linearized) fast variable $x$ at the Hopf point $y_H$ when no kick is delivered.  This will give us a reference unkicked, or ``natural" distribution $p_n(x)$ important in computing the buffer point below.  The trajectories of interest jump down from spiking when $y=y_i=y_{SN}$ and take $g^{-1}(y_H-y_{SN})/\e$ time units to reach $y_H$.  We use Eqn.~\eqref{e.gen} to write the resulting distribution as
\beq
p_{n}(x)=\tilde p(x(g^{-1}(y_H-y_{SN})/\e) | y_{SN}, (0,0)^T, \Sigma_{SN} )
\label{ss}
\eeq
where
\beqq
\Sigma_{SN}=\left(\begin{array}{cc} r_P(y_{SN})^2 & 0 \\0 & r_P(y_{SN})^2\end{array}\right)
\eeqq
and we have used Eqn.~\eqref{slow} to substitute in for time in Eqn.~\eqref{ss}.  The natural distribution is therefore defined to have an initial variance equal to the squared radius of the periodic orbits on $P$ when they vanish at $y_{SN}$, marking the end of the spiking phase.  Due to the long timescale of the slow dynamics, however, we find that the choice of initial variance has little effect on $p_{n}$.  We note that $p_{n}(x)$  is centered at $x=(0,0)^T$ with covariance depending on $\n$.

\begin{figure}[h]
\begin{center}
\includegraphics{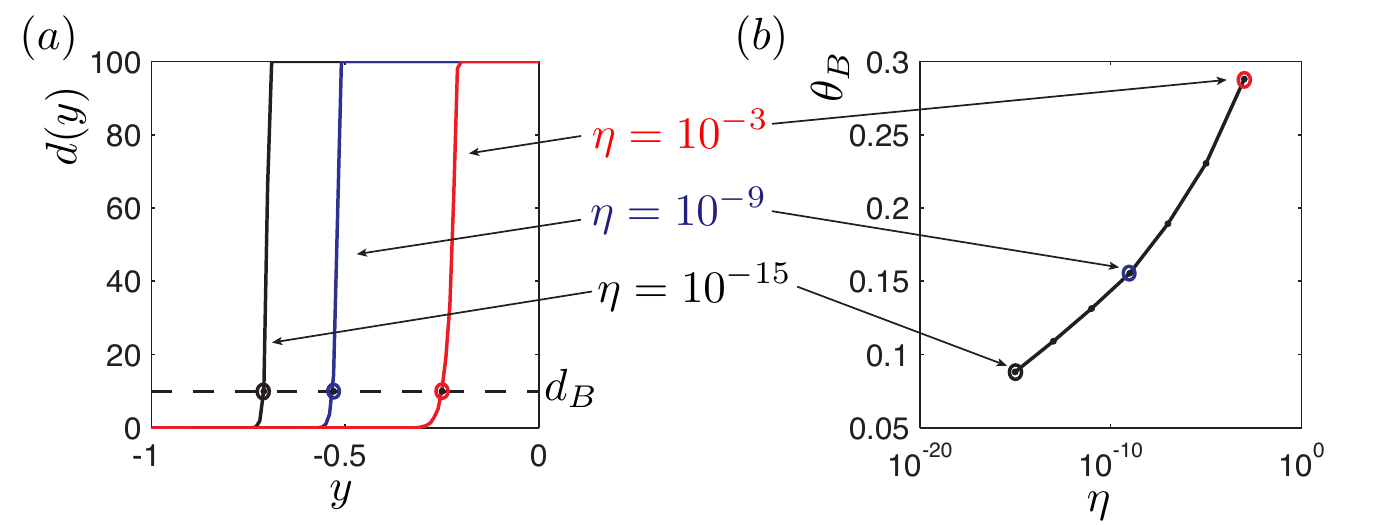}
\caption{(a) Plot of $d(y)$ with respect to $y$ for $\n=10^{-15},10^{-9},10^{-3}$. Dashed black line marks the threshold $d_B$ determining the points $y_B$ at the intersections with the $d(y)$ curve. (b) Plot of buffer phase point $\o_B$ with respect to noise strength $\n$, marking the onset of expansion in the weak kick map. Kick strength is $A=0.5$, as in Fig. \ref{fig:noisy_stats}.}
\label{fig:buffer}
\end{center}
\end{figure}

Next, we ask whether a trajectory that has received a weak kick can be expected to undergo comparable slow passage through the Hopf point as for the unkicked trajectories described by $p_{n}(x)$. Note here that an unkicked trajectory admits a maximal delay going through the Hopf bifurcation, and that some kicked trajectories will also have the same slow passage, as we expect $y_{SN}<y_B$.

If kicked trajectories typically pass through $y_H$ at locations $x$ with high probability density according to $p_{n}(x)$, we expect that they will have a comparable slow passage times as unkicked trajectories. On the other hand, if these trajectories are typically found where $p_{n}(x)$ is low, this indicates that they are further away from the branch of equilibria $S$.  Thus, they will tend to escape $S$ (i.e., jump up) sooner than the unkicked solutions. We obtain an approximation for the buffer point $y_B$ by asking when this distinction between kicked and unkicked trajectories at $y_H$ occurs.  

To this end, we compute the ``kicked" distributions of trajectories $p_{k,y}(x)$ for which a weak kick of amplitude $A$ is applied when $y(t)=y$.  Upon their arrival at $y_H$, we approximate these distributions as
\beq
  p_{k,y}(x) = \tilde p(x(g^{-1}(y_H-y)/\e) | y, (A,0)^T, \Sigma_{n} )
\label{dart}
\eeq
where the mean trajectory is at $(A,0)^T$ following the kick and $\Sigma_{n}$ is the covariance matrix from $p_n(x)$, under the assumption that the trajectory followed the natural burst cycle before the translation induced by the kick. We now assess to what extent $p_{k,y}$ and the natural distribution $p_{n}$ overlap. We use the symmetrized Kullback-Leibler Divergence between the two distributions 
\beq
d(y)=\frac{1}{2} \left(D_{KL}[p_n(x) \|p_{k,y}(x) ]+D_{KL}[p_{k,y}(x) \|p_n(x)] \right)
\label{dkl} 
\eeq
where
\beq
D_{KL}[p\|q]=\int_{\R^2}p(x)\ln\frac{p(x)}{q(x)}dx \;.
\eeq

For $y$ sufficiently far from $y_H$, the distribution $p_{k,y}$ has enough time to converge close to $p_n(x)$ before $y(t)$ reaches $y_H$.  As a consequence, $d(y)$ will be close to zero. We define the $A$-dependent buffer point $y_B$ to be the first value of $y$ for which $d(y)$ grows beyond a threshold $d_B$. Since the distributions $p_{k,y}$ and $p_n$ are sharply peaked Gaussians (with variance of order $\eta^2$), $d(y)$ quickly explodes -- to several orders of magnitude above one -- when the two distributions fail to overlap. Therefore, we choose $d_B=10^1$ as a good indication of separation among kicked vs. unkicked trajectories, as illustrated in Fig. \ref{fig:buffer} (a).

Finally, we determine the phase point $\o_B$ corresponding to $y_B$, via \eqref{phase_def}:
\beq
\o_B=\frac{h_S^{-1}(y_B)}{T}
\label{onset}
\eeq
where $T$ is the mean period of the unperturbed trajectory.  This point marks the onset of expansion for the associated kick map.  We find excellent agreement of this prediction with numerically computed kick maps, as seen from Figs. \ref{fig:buffer} (b) and \ref{fig:noisy_stats} (b).

\subsubsection{Deriving the kick map for noisy bursters}


We next derive an approximate expression for the complete kick map in the presence of weak noise.  Our first step is to decouple the $y$ dynamics from the (noisy) fast variables.  This is guided by the assumption that, due to weak noise, most trajectories closely follow $S$ and $P$. We then proceed to derive $h_S$, $h_S^{-1}$ and $h_P^{-1}$ as done in Sect. \ref{map}, with a single modified value:  $y_J$. Indeed, as described above, the jump up point $y_J$ is closer to $y_H$ for noisy bursters and we numerically compute its value for each noise strength $\n$.

In the presence of a buffer point, a weak kick can now have two outcomes:  either it has no effect if it is received when $0<y<y_B$ as it does not alter slow passage, or it shortens slow passage as in the deterministic case, when $y_B<y<y_w$. To capture the slow passage effects of weak kicks (responsible for the expansion in our map) in the presence of noise, we still use integral condition~\eqref{integral_condition}. Although this formula was derived for deterministic bursters, it relies in the linearization of the fast dynamics about $S$, which we assume remains valid even in the noisy case. As a result, numerical simulations show that $\tilde{y}_j(y)$ holds true for $y \in (y_B,y_H)$ except in a short interval to the right of $y_B$ where small errors are observed. Note that these errors diminish with smaller noise for which $\tilde{y}_j(y_B)$ is quite close to the numerically computed $y_J$. We proceed to write an expression for our new kick map, which now contains an addition piecewise-defined section arising from the presence of 
 $\o_B$: 
\begin{equation}
F_{A}(\o) = \left\{ \begin{array}{ll}
\o & \textrm{if $\o \in [0,\o_B]$}\\
\o+h_P^{-1}(\tilde{y}_j(h_S(\o)))-h_S^{-1}(\tilde{y}_j(h_S(\o))) & \textrm{if $\o \in [\o_B,\o_w)$}\\
h_P^{-1}\circ h_S(\o) & \textrm{if $\o \in [\o_w,T_S]$}\\ 
\o & \textrm{if $\o \in [T_S,1]$}.
\end{array} \right.
\label{noisy_map}
\end{equation}
We obtain excellent fits as shown in Fig. \ref{fig:noisy_stats} (b). We end by noting that we get similar fits for various kick strengths $A$ as well as distinct parameters sets (i.e. $b>0$) for the normal form model (not shown).

\subsection{Effect of noise on iterated dynamics}
\label{noisy_it}

\begin{figure}[h]
\begin{center}
\includegraphics{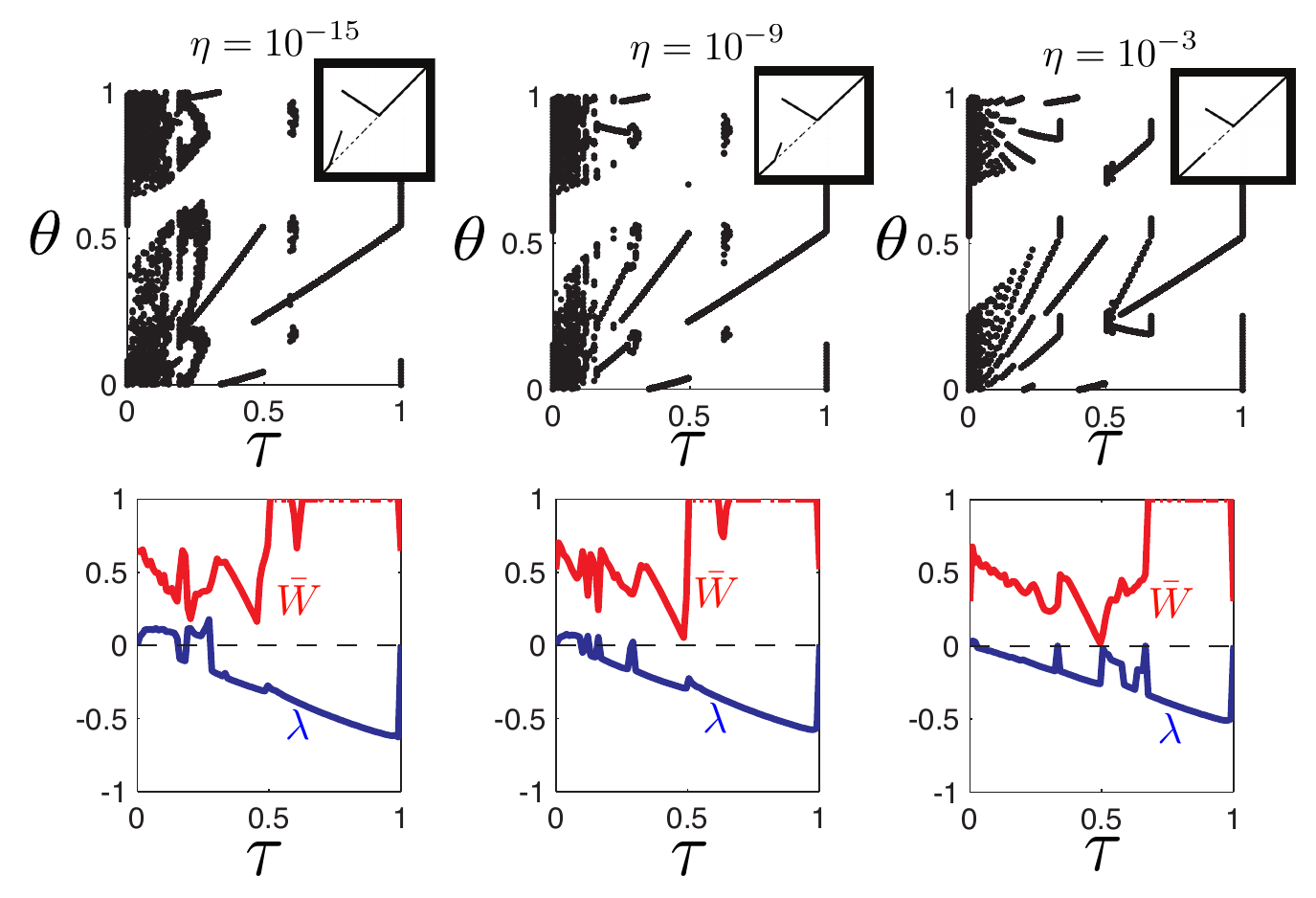}
\caption{Orbit diagrams and related measures for weak kick maps with $A=0.5$, for various noise strengths $\n$ (see Fig. \ref{fig:noisy_stats}). Top to bottom:  orbit diagrams of 100 cells after 150 iterations;  synchrony measure $\bar{W}$ (red) and averaged Lyapunov exponent $\l$ (blue). Left to right:  low noise strength $\n=10^{-15}$;  medium noise strength $\n=10^{-9}$;  high noise strength $\n=10^{-3}$.}
\label{fig:noisy_bif}
\end{center}
\end{figure}

We now explore the dynamical properties of the kick maps computed in the presence of noise. Figure \ref{fig:noisy_bif} shows orbit diagrams, synchrony measures and averaged Lyapunov exponents three maps (computed as for Fig.~\ref{fig:bif_diag}). For $\n=10^{-5}, 10^{-9}$, the maps retain some expansion and we see behavior that appears chaotic for small positive values of $\tau$. This range of $\tau$ values shrinks as the as the expansive region of the map gives way to neutrality  with increasing noise; at the same time, the ``support" of the orbit diagrams appears to decrease.
Arguments similar Lemma \ref{lemma_chaos} can be formulated for these cases to show the existence of positive Lyapunov exponents (see also bottom panels of Fig.~\ref{fig:bif_diag}). 

For the map with $\n=10^{-3}$, there is no expansive region:  although there is still a slow passage effect the kick's amplitude $A=0.5$ induces a cutoff $\o_w$ small enough such that the system does not retain memory of any kick ($\o_w<\o_B$).  Thus, we cannot expect positive Lyapunov exponents (see bottom panel in Fig.~\ref{fig:noisy_bif}). However, the orbit diagram shows a broad spread of points for small positive $\tau$. These are stable, high period orbits, originating from border crossing bifurcations as $\tau$ increases. This has been established by B\'elair in the context of periodically forced integrate and fire oscillators \cite{Belair:1986p8707}, where a very similar map is studied:  he shows that as the map is shifted vertically, stable periodic orbits of  a wide range of periods can be found, following a Farey tree sequence.

In sum, small positive values of $\tau$ result in either positive Lyapunov exponents or high period orbits for the weak kick maps, depending on the noise level $\eta$ assumed in deriving the (deterministic) map.  In the first case, expansion directly desynchronizes cells; in the second, we will see that the high period of orbits, coupled with additional variability due to the underlying noise, can have a similar effect.

We now introduce stochastic terms into our discrete kick map dynamics to account for the variability in burst periods discussed above (i.e., $CV \neq 0$). 
If, for a given cycle, a cell has a shorter/longer period than the one used to compute it's kick map, its phase following a kick will be slightly shifted from the phase given by the iteration of the map. To capture this, we introduce \textit{jitters}:  additive stochastic terms acting on $\tau$, independent for every cell. The goal is not to capture the exact phase response of cells, but rather to give a qualitative account for the impact of period variability on statistical metrics such as our synchrony measure. We proceed as follows:  

Since the construction of our kick maps rescales the period of any burster to be 1, the $CV$ can be interpreted as the standard deviation of a burst cycle's period. We define a \textit{jitter} $\zeta$ to be random variable drawn from a normal distribution with zero mean and standard deviation equal to the $CV$ of the case we are considering. Suppose we want to model the phase evolution of $M$ cells subject to a common, periodic kick train of period $\tau (\rm{mod}~ 1)$, of amplitude A, subject to stochastic forcing of strength $\n$. The phase of the $m^{th}$ cell right before the $n+1^{st}$ kick can be written as follows
\begin{equation}
\o_{n+1}^m=F_{A,\n,(\tau+\zeta^m_n)}(\o_{n}^m)
\label{noisy_iteration}
\end{equation}
where $\zeta_n^m \sim^{iid}N(0,CV_\n)$.  In other words, at every iteration, we draw a different jitter $\zeta$ for every cell. Note that we modified our notation to emphasize the map's dependence on noise strength $\eta$ (Fig.~\ref{fig:noisy_stats}). 

\begin{figure}[h]
\begin{center}
\includegraphics{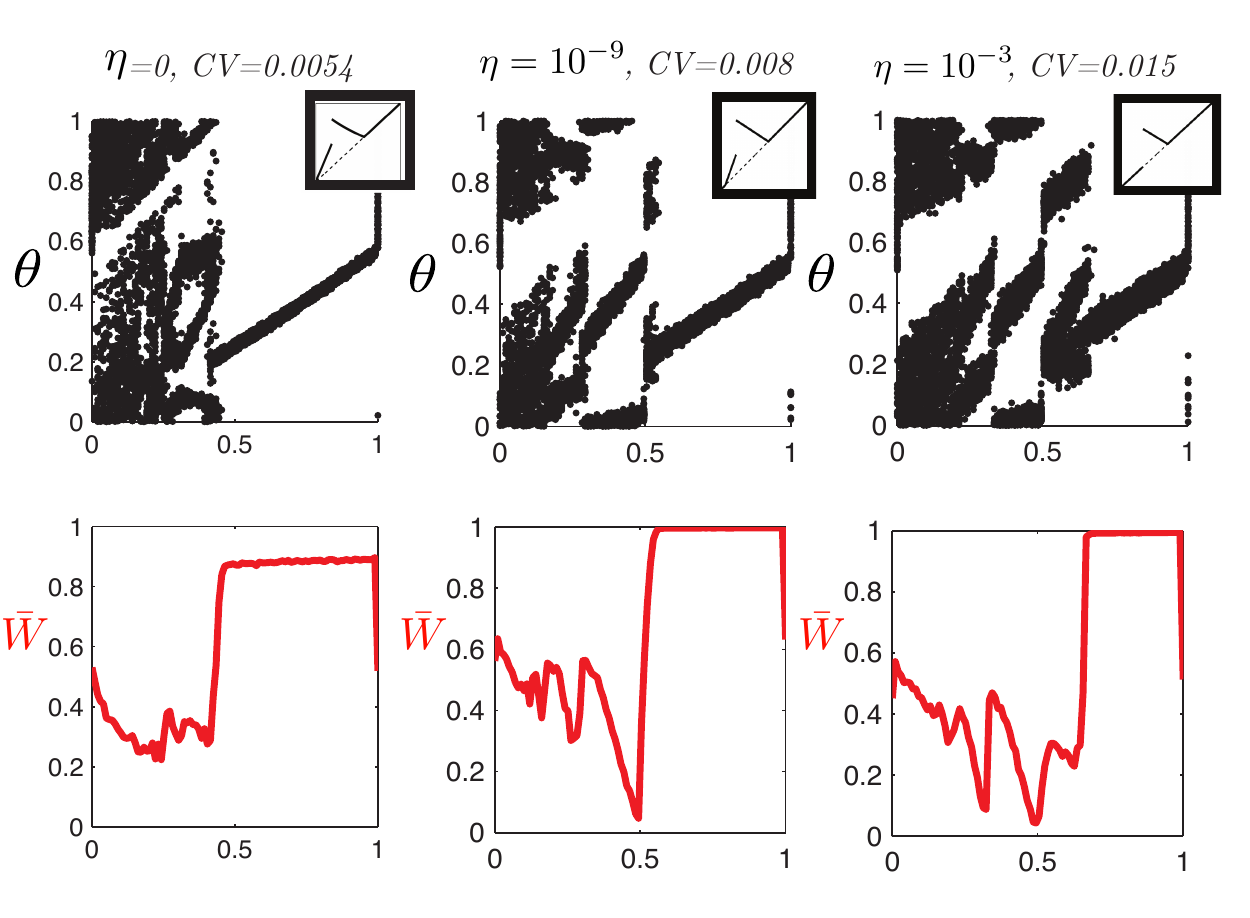}
\caption{Orbit diagrams and synchrony measure for weak kick maps with $A=0.5$, for various noise strengths $\n$ and added jitters $\zeta(CV)$. Top to bottom:  orbit diagrams of 100 cells after 150 iterations;  synchrony measure $\bar{W}$. Left to right:  low noise strength $\n=10^{-15}$;  medium noise strength $\n=10^{-9}$;  high noise strength $\n=10^{-3}$.}
\label{fig:jit_bif}
\end{center}
\end{figure}

Using (\ref{noisy_iteration}), we again iterate 100 cells 150 times with added jitters and plot the orbit diagrams and synchrony measures in Fig. \ref{fig:jit_bif}, for the same three levels of $\eta$ as in the preceding figures. 
Jitters, as expected, ``smear" orbit diagrams, with a greater effect for larger $\eta$.  In particular, note the smoothing of periodic points for small positive $\tau$ in the high noise case ($\n=10^{-3}$).  Interestingly, the twin effects of noise in reducing expansion but increasing cell-to-cell jitter result in comparable levels of the synchrony measure $\bar W$ across the three cases.

\bigskip

We reiterate our main conclusion:  although the underlying mechanisms differ across a wide range of noise strengths $\eta$, pulsatile inputs in the ``weak" kick regime -- with an input frequency slightly slower than (a multiple of) cells' intrinsic frequencies -- will result in desynchrony among a population of recipient cells.

\section{Validity of phase reduction, O.D.E. simulations, and a neurobiological model}
\label{numerics} 

In this section, we explore the validity of our phase reductions in both deterministic and noisy cases -- and our analysis of the discrete kick map that follows -- by numerically integrating the O.D.E.s themselves.  Rather than demonstrating a complete correspondence between the kick map and solutions of the differential equations, we seek to verify that an informed choice of kick amplitude and period, based on the kick maps, does indeed yield the predicted (de)synchrony behavior among solutions to the O.D.E.s.   Specifically, we show that small, positive values of $\tau$ (i.e., $\tau \in (0,\tau_C)$ in the deterministic case) with amplitudes in the ``weak" regime lead to the greatest desynchrony; conversely, large values of $\tau$ synchronize cells.

We first consider the normal form system (\ref{canon_bautin}), and then turn to a biologically detailed neuronal model for which our main findings persist.  All numerical computations were carried out in MATLAB. We use the stiff solver \verb"ode15s" with both absolute and relative tolerances set to $10^{-6}$ to integrate all differential equations; input kicks and additive noise are treated as non-autonomous terms by the solver.

\subsection{Normal form model}
\label{canon_numerics}

Here, we numerically integrate a population of $N=30$ uncoupled cells governed by system (\ref{canon_bautin}) or its stochastic counterpart Eqn.~(\ref{canon_noisy}), taking the ``large" noise value $\n=10^{-3}$ studied above (with \{$\e=0.01$, $w=1$, $a=0.8$, $b=0$\}). We concentrate on one weak kick amplitude, $A=0.5$. In each case, we implement periodic kicks that correspond to $\tau=0.1$ and $\tau=0.8$ for the kick map, to illustrate desynchronizing and synchronizing behavior respectively.  More precisely, we use a kick period equal to $(1+\tau) \times T$ where $T$ is the natural period of the O.D.E.'s  burst cycle.  This allows trajectories at least one natural period to relax toward the unperturbed cycle in between kicks (similar dynamics occur for periods $T\times (n+\tau)$, $n \in \N$). For $\n=0$, $T \simeq 465$ while for $\n=10^{-3}$, $T \simeq 337$ (Fig.~\ref{fig:noisy_stats}(d)). 

\begin{figure}[h]
\begin{center}
\includegraphics{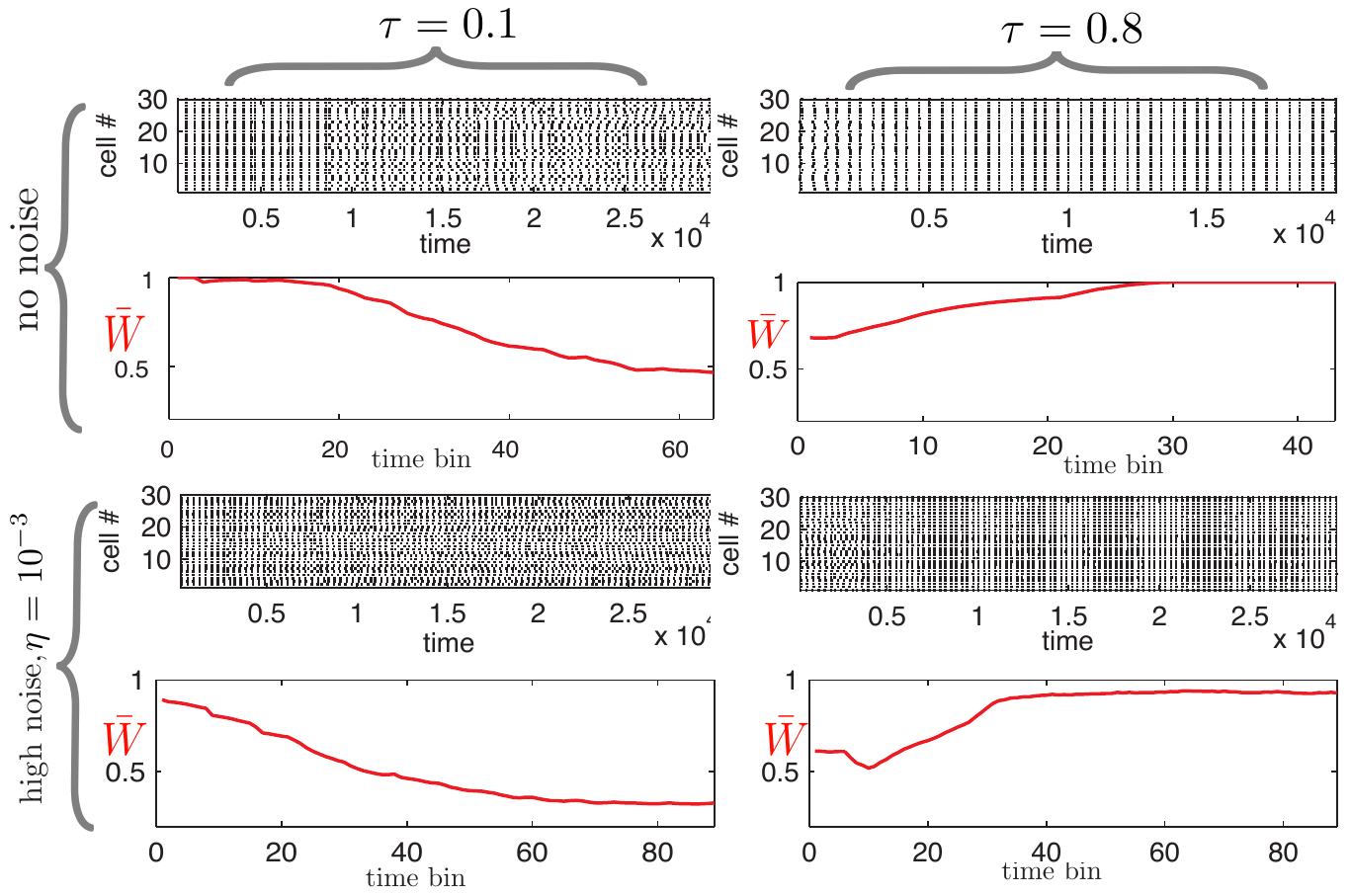}
\caption{A population of 30 numerically integrated solutions of Eqn.~(\ref{canon_bautin}) or (\ref{canon_noisy}), all receiving a common, periodic weak kick input ($A=0.5$). Left column: kick period $T \times 1.1$ (equivalent to $\tau=0.1$) results in population desynchrony. Right column:  kick period  $T \times 1.8$ (equivalent to $\tau=0.8$) synchronizes the population. Top row:  no noise ($\n=0$). Bottom row:  high noise ($\n=10^{-3}$). Black dots give the raster plot (see text); red curves plot the synchrony measure $\bar{W}$ vs. time (see text).}
\label{fig:no_noise_stim}
\end{center}
\end{figure}

Figure~\ref{fig:no_noise_stim} displays the results via {\it raster plots}:  for each cell, a dot is placed at the moment that spiking terminates (corresponding to phase $\theta=0$).  We also plot a synchrony measure for the simulated population.  This is done by assigning phases to each cell relative to their most recent spike termination event, as a fraction of elapsed time partitioned in bins of length $T$.
To better illustrate the desynchronizing effect of weak kicks with $\tau=0.1$, initial conditions are chosen at random with phases at most $2\%$ apart (i.e., an initially synchronized population); to illustrate the synchronizing effect of kicks with $\tau=0.8$, initial phases are allowed to be more sparse. In all cases, we let the cells evolve without inputs for a few burst cycles, and then begin to apply the pulsatile inputs.

The results agree well with predictions from the kick map:  a weak kick of $A=0.5$ administered at $T \times 1.1$ successfully spreads cells apart while the same kick with period $T \times 1.8$ synchronizes the population.  Note that we chose these values of $\tau$ only as informed guesses; other nearby values can achieve similar results in both synchronizing and desynchronizing the population. Moreover, results for various other kick amplitudes also agree well with the behavior predicted from the associated kick maps (not shown).

\subsection{GPe bursting neuron}
\label{GPe}
 We now investigate whether the mechanisms described above will persist for a more biologically detailed model.  Specifically, we study a 5-dimensional, Hodgkin-Huxley-type model of a neuron from the GPe basal ganglia nucleus~\cite{Terman:2002p125,Best:2007p128}.  This model produces elliptic bursting where the onset of spiking is due to a subcritical Hopf bifurcation and a burst termination is due to a saddle node on an invariant circle, in agreement with the normal form system (\ref{canon_bautin}). In detail, the fast variables are the voltage $V$, potassium current gating variable $n$, sodium current gating variable $h$, and calcium T-current gating variable $r$. The slow variable is calcium concentration $Ca$. The equations are as follows:
 \begin{equation}
\begin{split}
  C_m \frac{d V}{d t} & =-I_{Ca}+I_{Na}+I_K+I_L+I_{AHP}+I_T-I_{app}+I(t)+\n\xi(t)\\
  \frac{d n}{d t} & =-\frac{  \phi_n{(n-n_\infty)}}{\tau_n}\\
  \frac{d h}{d t} & =-\frac{ \phi_h {(h-h_\infty)}}{\tau_h}\\
  \frac{d r}{d t} & =-\frac{{\phi_r(r-r_{\infty})}}{\tau_r}\\
  \frac{d Ca}{d t} & =-\e {(I_{Ca}+I_T+ k_{Ca} Ca)} \\
\end{split}
\label{eqn:Terman}
\end{equation}
where the $I$ terms represent membrane currents and are functions of the gating variables and the voltage; all definitions and parameter values are as in \cite{Terman:2002p125}.  Additionally, the terms $I(t)$ and $\n\xi(t)$ represent the pulsatile inputs and the noise term, entering as currents. 
\begin{figure}[h]
\begin{center}
\includegraphics[width=400pt]{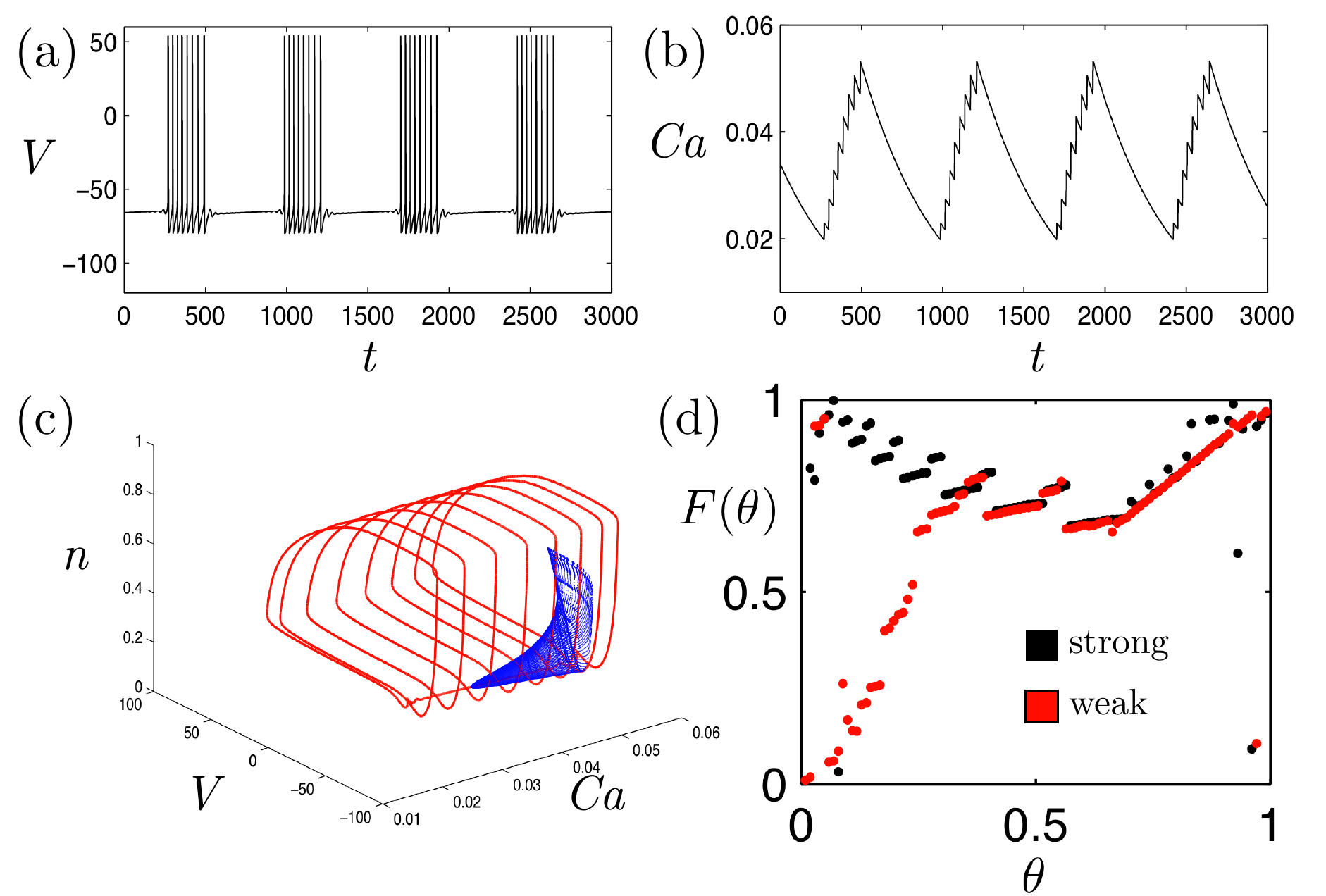}
\caption{Plots from model (\ref{eqn:Terman}). (a) Voltage trace. (b) Calcium trace. (c) Bursting solution in $n$, $V$, $Ca$ plane in red and separatrix $U$ in blue. (d) Two numerically computed kick maps for model (\ref{eqn:Terman}) with high noise strength $\n=10^{-3}$;  strong kick ($A=30$) in black and weak kick ($A=3$) in red.}
\label{fig:Terman}
\end{center}
\end{figure}

Figure \ref{fig:Terman}(a) shows that the shape of action potentials and the timescales differ from those of the normal form model (\ref{canon_bautin}); nevertheless, the dynamics have a very similar structure.  In particular, panel (c) shows a (projected) 3-dimensional plot of a bursting trajectory together with the skewed separatrix $U$, computed using the MATCONT package \cite{Dhooge:2003p9696}.  Figure \ref{fig:Terman}(d) shows two numerically computed kick maps for both strong and weak kick amplitudes ($A=30$ and $A=3$, respectively). These maps were computed in the presence of noise with $\n=10^{-3}$, which is the largest noise strength that keeps $CV$ at $O(10^{-2})$ for the simulations; maps represent the average phase response taken over ten runs with different noise realizations.

Overall, the structure of these maps is more complex than for the normal form model. In particular, ``small" plateaus and associated  discontinuities are promienent. As for the normal form model, there are as many plateaus as there are spikes in an unperturbed burst, and kicked solutions that elicit a certain number of spikes in the subsequent burst accumulate in each plateau. However, for this model, the slow variable (calcium concentration) varies more during a spike and creates bigger gaps between plateaus. As a result, even for a strong kick, certain values of $\tau$ yield localized stable periodic orbits, as opposed to only fixed points. These appear via border collision bifurcations due to discontinuities between plateaus (not shown). However, the small amplitude of these periodic orbits keeps the cells attracted to them quite synchronized. 

Additionally, the shape of the left part of the weak kick map is also quite distorted compared with maps derived from the normal form model.  In particular, notice that there are large discontinuities close to zero. This is due to the skewed cone shape of the separatrix $U$:  since the neuron model does not have the same symmetry as our normal form system, when the solution drops down from spiking, it spirals towards the resting branch and some lobes of this spiral come very close to the separatrix. When the solution is kicked, even weakly, on the upper part of a lobe, it passes the separatrix and jumps to the spiking state; the same weak kick will not have this effect if it is delivered only moments later.  The resulting large gaps in the weak kick map add to the complexity of the dynamics for low $\tau$.

\begin{figure}[t]
\begin{center}
\includegraphics{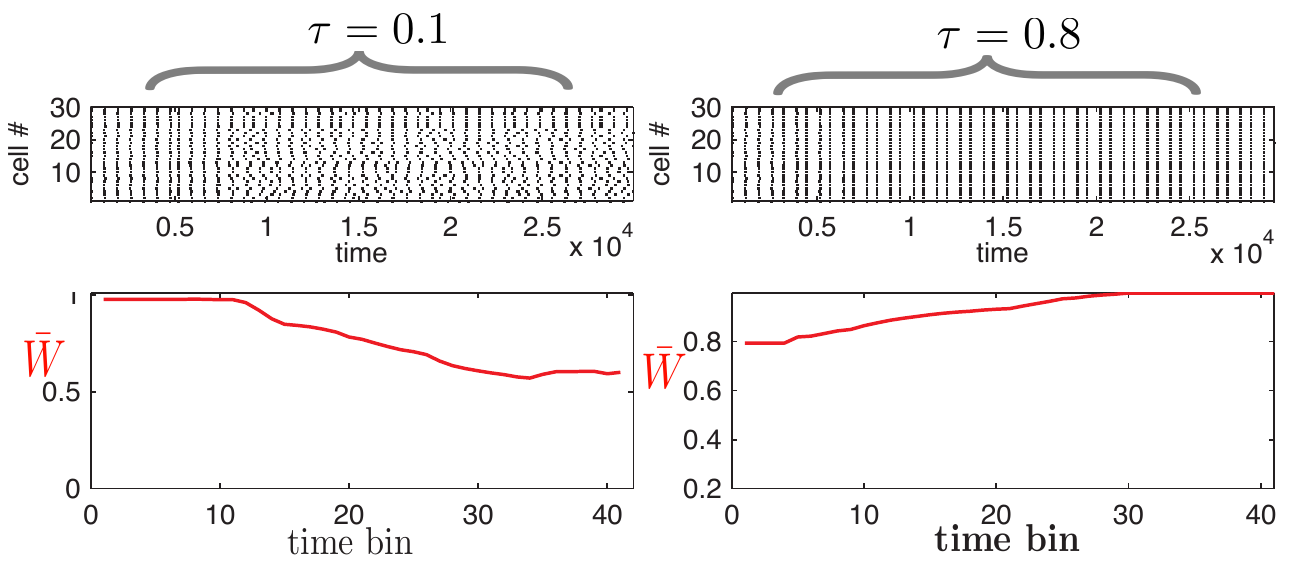}
\caption{A population of 30 numerically integrated solutions of Eqn.~(\ref{eqn:Terman}) with noise strength $\n=10^{-3}$, receiving a common periodic weak kick input ($A=3$).  Left column: $T \times 1.1$ (equivalent to $\tau=0.1$) results in population desynchrony. Right column:  kick period  $T \times 1.8$ (equivalent to $\tau=0.8$) synchronizes the population. Other plotting details also as for Fig.~\ref{fig:no_noise_stim}.}
\label{fig:terman_sim}
\end{center}
\end{figure}

Apart from these differences, the prominent features observed in the kick map for the normal form model remain:  the neutral/expanding left branch for the weak kick map, and the contracting middle branch and neutral right branch. We repeat the numerical experiment described in Sect. \ref{canon_numerics}, this time only for the noisy case ($\n=10^{-3}$), and plot the results in Fig. \ref{fig:terman_sim}.  We see the expected synchronization and desynchronization from weak periodic kicks ($A=3$) with periods equivalent to $\tau=0.8$ and $\tau=0.1$ respectively. 
While cells do not appear to become as desynchronized for the $\tau=0.1$ case as in the normal form model, it is reasonable to believe that a more detailed analysis of the kick map for the neural model could identify ($A,\tau$) combinations that would further desynchronize cells.

\section{Composition of multiple periodic inputs, and an application to DBS}
\label{composition}

Above, we showed how weak, periodic inputs can lead to desynchronization for populations of uncoupled bursting cells.  But how well can such inputs compete with other, synchronizing effects?  The answer is important in varied applications.  A prominent one is Deep Brain Stimulation (DBS) therapy for Parkinson's disease.  Here, pathologically high levels of synchrony occur among bursting cells in the basal ganglia.  Synchrony in some basal ganglia areas is in large part driven by common, periodic inputs from other areas (see~\cite{Rubin:2004p135} and references therein).  
A DBS electrode delivers pulsatile electrical signals that are designed to mitigate the effects of this synchrony.  Thus, {\it two} common periodic inputs are received by bursting neurons, possibly with competing effect.

Using the normal form model \eqref{canon_bautin}, we undertake a brief demonstration of how our results could be applied to this setting, for the GPe basal ganglia nucleus that contains neurons believed to be elliptically bursting.  We do not attempt detailed, biologically complete modeling or aim for direct clinical relevance, and as such note several limitations.  First, the source of intrinsic entrainment here is a purely common, periodic drive;  lateral connections between bursting cells, believed to be sparse and weak in Parkinsonian regime~\cite{Rubin:2004p135,Pirini:2009p521}, are neglected.  Second, GPe is not the most common target for DBS in practice, though it has been the focus of several emerging studies \cite{Johnson:2009p1400,Pirini:2009p521}.  Nevertheless, a better theoretical understanding of the interactions between intrinsic and applied inputs to the GPe could, in the long term, contribute to the design of signals that desynchronize {bursting} neurons by targeting key instabilities~(cf.~\cite{Hauptmann:2007p2488,Hauptmann:2009p2319,Danzl:2009p7990,Feng:2007p148}), taking inspiration from similar findings for oscillatory neurons~\cite{Danzl:2009p7990,Goldobin:2005p9724,Kosmidis:2003p9761,Winfree:2001p9818,Lin:2006p9914}.

We suppose that a population of bursters receives a first sequence of synchronizing periodic impulses with period $\tau_1$ and ``strong" amplitude $A_1$.  The action of these inputs on burst phases is given by the kick map $F_{A_1,\tau_1}(\o)$.  As throughout our paper, this returns the phase of a cell following a kick, $\tau_1$ time units later.  Aiming to counteract the synchrony due to the first kick sequence, we introduce a second series of kicks of strength $A_2$.  We assume that these have the same period, but are delayed by an amount $\tau_2$.  That is, the cell receives a $A_2$-kick $\tau_2$ time units following each $A_1$-kick. We wish to write the kick map that captures the effect of such doublets of kicks. 

In this context, the shift-time following a $A_1$-kick must be taken to represent the phase of cells right before the $A_2$-kick and the first application of the map must be $F_{A_1,\tau_2}(\o)$. Similarly, we must shift the $A_2$ map by $\tau_1-\tau_2$ to retrieve the phase of a cell before the next $A_1$-kick. Note that neither $\tau_2$ nor $\tau_1-\tau_2$ should be too small for this map to be valid, specifically in the presence of weak kicks when the cell must have time to enter its spiking phase before the following kick, for the map we derive to remain valid.  When this restriction is satisfied, the doublet map is given by
\beq
F_{A_1A_2,\tau_1\tau_2}= F_{A_2,\tau_1-\tau_2}\circ F_{A_1,\tau_2} \;.
\label{doublet_map}
\eeq


\begin{figure}[ht]
\begin{center}
\includegraphics{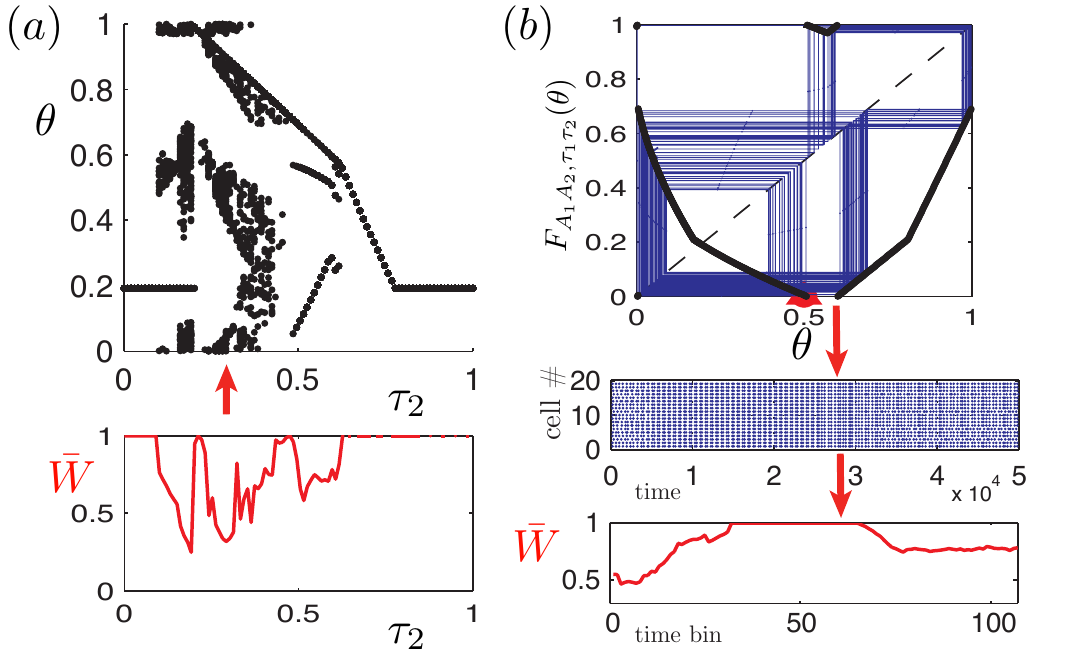}
\caption{(a) Orbit diagram and synchrony measure for $F_{A_1 A_2,\tau_1 \tau_2}$ while varying $\tau_2$. (b) Top : cobweb diagram of $F_{A_1 A_2,\tau_1 \tau_2}$ with $\tau_2$ indicated by the red arrow in (a). (b) Bottom : O.D.E. simulation of 20 cells initially synchronized by a strong input and then desynchronized by competing weak kicks (starting at red arrows).}
\label{fig:composed}
\end{center}
\end{figure}

An example is shown in Fig.~\ref{fig:composed}.  Suppose we start with an entraining input of strong kicks with $A_1=1.5$ and $\tau_1=0.4$. We seek to oppose this synchronizing effect with weak kicks of amplitude $A_2=0.5$. We use the two first maps of panel (a) of Fig. \ref{fig:computed_map} to build the resulting doublet map $F_{A_1A_2,\tau_1\tau_2}$. In panel (a) of Fig.~\ref{fig:composed}, we compute the orbit diagram of this map (as done in Sect. \ref{dynamics}) while treating $\tau_2$ as our variable parameter. Using this diagram, we select $\tau_2=.375$ (marked by a red arrow), associated with a low synchrony measure.  We plot a cobweb diagram of a sample orbit in the top of panel (b) of Fig.~\ref{fig:composed}, which clearly demonstrates the destabilizing effects of expansive regions in the doublet map.  

We then verify the properties desynchronization predicted by the doublet map by numerically solving the underlying O.D.E.s \eqref{canon_bautin}, for twenty model cells.    We begin with initial conditions such that the phases are desynchronized and apply strong kicks ($A_1=1.5$) at $1.4 \times T$, where T is the natural period of the bursters (i.e., corresponding to  $\tau_1=0.4$). As predicted, the cells synchronize in response; see the binned synchrony measure rising up to one in the bottom of panel (b) of Fig.~\ref{fig:composed}, or the raster plots above. After synchrony has developed (red arrows in panel (b)), we ``switch on" the sequence of weak kicks, leaving in place the original strong kick sequence.  Weak kicks are applied $0.375 \times T$ time units following each strong kick (i.e.,  $\tau_2=0.375$).  The desynchronizing impact predicted by the doublet map is clear in both the scatter in raster plots and in the drop in the synchrony metric $\bar W$ that develops after the weak input begins to be applied.


\section{Conclusion}
\label{discussion}

We study the behavior of a population of identical elliptic bursters receiving a periodic sequence of pulsatile inputs, or kicks.  Our aim is to understand which input sequences will result in desynchronized vs synchronized bursts across the population.
Following and extending the approach in \cite{Best:2007p128}, we first conduct a phase reduction of the burst dynamics to a circle map, using a slow/fast decomposition.  This ``kick map" depends on two parameters -- the kicks' amplitude $A$ and (relative) period $\tau$ -- and maps phases from their states just before one input pulse to their states just before the next pulse arrives.  We next study the effect of varying $A$ and $\tau$ using a normal form model for elliptic bursting (Eqn.~(\ref{canon_bautin})). 

We find that for strong kicks -- i.e., with $A$ sufficiently large so that the cell will always be spiking following an input -- almost any choice of kick period $\tau$ resulted in $1:1$ phase locking, and hence synchrony across the population. For weaker kicks, we find a rich dynamical structure. In particular, the interaction of a weak perturbation with the slow passage effect through a subcritical Hopf point induces an expansive region in the kick map. By varying the kick period, we witness the appearance of stable fixed points, periodic orbits and regimes with positive Lyapunov exponent.  As expected, this leads to desynchronization of the population.  Overall, we divide the ($A$, $\tau$) parameter space into the three regions shown in Fig.~\ref{fig:trajectories}(a), corresponding to unstable, desynchronizing dynamics, 1:1 phase locking, and intermediate, complex behavior,  The former, desynchronizing regime is associated with relatively weak kicks of periods slightly slower than the natural burst period ($0<\tau<\tau_C$, see Eqn.~\eqref{tau_c}). 

We also study the effect of stochastic perturbation via noise terms. We find that the phase reduction retains its validity but the kick map changes shape, presenting less expansion as the noise increased.  Importantly, population desynchrony still results from weak kicks with comparable values of $\tau$ in this case, but through a different mechanism than the instabilities that occur for the noise free case.  Here, desynchrony follows from a combination of high-period orbits and the noise itself.  Overall, this phenomenon is related to the discontinuous nature of the circle map at hand; 1:1 phase locking rather than the complex dynamics observed would be expected for small $\tau$ for many smooth maps~\cite{Glass:1983p9985,Glass:1985p10448}.  
 
 We then test the predictions of the reduced circle maps via numerical simulation of the original O.D.E. system, finding qualitative agreement.  Additionally, we simulate a more biologically realistic model of a GPe neuron, and continue to find agreement with the general predictions of our maps.
Finally, we show that it is possible to use the kick map framework to study the effect of multiple sequences of inputs to a cell population. We build an example showing that carefully timed weak kicks can compete with an entraining strong input to successfully desynchronize a population of bursting cells.

As a closing remark, we note that the kick map studied here can also capture the effect of  pulsatile input signals that are neither periodic, nor have a fixed kick amplitude. For any given sequence $\{A_n,\tau_m\}$, where $A_n$ is the amplitude of the $n^{th}$ kick and $\tau_n$ is the delay between kicks $n$ and $n+1$, the relevant system is the composition of the maps $F_{A_n,\tau_n}(\o)$. This gives rise to an 
iterated function system (IFS) acting on $S^1$. There is a growing body of literature dealing with these objects and their application to this problem could eventually help us to understand the behavior of bursting cells under arbitrary -- and possibly stochastic -- stimulation patterns.

\section{Acknowledgements}  

We thank Jonathan Rubin, David Terman, Pablo G. Barrientos, Artem Raibekas, and Joshua Goldwyn for helpful discussions and insights.  This research was supported by NSF grant DMS-0817649, by a NSERC Graduate Fellowship (to G.L.) and by a Career Award at the Scientific Interface from the Burroughs-Wellcome Fund (to E.S.-B.). 
  
   
   \begin{small}
   \bibliographystyle{plain}
   \bibliography{maBibli}
   \end{small}

\end{document}